\documentclass[12pt]{amsart}
\usepackage{}
\usepackage{txfonts}
\usepackage{amssymb}
\usepackage{eucal}
\usepackage{graphicx}
\usepackage{amsmath}
\usepackage{amscd}
\usepackage[all]{xy}           %xypic macro for latex2.09
\usepackage{tikz}
\usepackage{amsfonts,latexsym}
\usepackage{xspace}
\usepackage{epsfig}
\usepackage{float}
\usepackage{color}
\usepackage{fancybox}
\usepackage{colordvi}
\usepackage{multicol}
\usepackage{colordvi}
\usepackage[shortlabels]{enumitem}
\usepackage{ifpdf}
\ifpdf
  \usepackage[colorlinks,final,backref=page,hyperindex]{hyperref}
\else
  \usepackage[colorlinks,final,backref=page,hyperindex,hypertex]{hyperref}
\fi
\usepackage{tikz}
\usepackage[active]{srcltx}

\makeatletter

\theoremstyle{definition}
\newtheorem{theorem}{Theorem}[section]
\newtheorem{prop}[theorem]{Proposition}
\newtheorem{defn}[theorem]{Definition}

\newtheorem{prop-def}{Proposition-Definition}[section]

\newtheorem{remark}[theorem]{Remark}

\newcommand{\nc}{\newcommand}

\nc{\delete}[1]{{}}
\nc{\mmargin}[1]{}

\delete{
\nc{\mlabel}[1]{\label{#1}}  % Use this to suppress names
\nc{\mcite}[1]{\cite{#1}}  % Use this to suppress names
\nc{\mref}[1]{\ref{#1}}  % Use this to suppress names
\nc{\mbibitem}[1]{\bibitem{#1}} % Use this to show number
}

\delete{
\nc{\mlabel}[1]{\label{#1}  % Use the next two lines to show names
{\hfill \hspace{1cm}{\bf{{\ }\hfill(#1)}}}}
\nc{\mcite}[1]{\cite{#1}{{\bf{{\ }(#1)}}}}  % Use this lines to show names
\nc{\mref}[1]{\ref{#1}{{\bf{{\ }(#1)}}}}  % Use this lines to show names
\nc{\mbibitem}[1]{\bibitem[\bf #1]{#1}} % Use this to show name
}
\delete{
\nc{\mlabel}[1]{\label{#1}}  % Use this to suppress names
\nc{\mcite}[2][]{\cite[#1]{#2}}  % Use this to suppress names
\nc{\mref}[1]{\ref{#1}}  % Use this to suppress names
\nc{\mbibitem}[1]{\bibitem{#1}} % Use this to show number
}

\nc{\mlabel}[1]{\label{#1}  % Use the next two lines to show names
{\hfill \hspace{1cm}{\bf{{\ }\hfill(#1)}}}}
\nc{\mcite}[1]{\cite{#1}{{\bf{{\ }(#1)}}}}  % Use this lines to show names
\nc{\mref}[1]{\ref{#1}{{\bf{{\ }(#1)}}}}  % Use this lines to show names
\nc{\mbibitem}[1]{\bibitem[\bf #1]{#1}} % Use this to show name
\nc{\difforg}{U_D}
\nc{\rbforg}{U_{RB}}
\nc{\multforg}{U_{M}}
\nc{\diffree}{{}_DF}
\nc{\difcof}{F_D}
\nc{\rbfree}{{}_{RB}F}
\nc{\rbcof}{F_{RB}}
\nc{\multfree}{{}_{M}F}
\nc{\multcof}{F_M}
\nc{\rc}{\prec}
\nc{\lc}{\succ}
\nc{\rcc}{\ll}
\nc{\lcc}{\gg}
\nc{\rrc}{\preceq}
\nc{\llc}{\succeq}
\nc{\rs}{\acute{\succ}}
\nc{\ls}{\grave{\succ}}
\nc{\lj}{\grave{\prec}}
\nc{\rj}{\acute{\prec}}
\nc{\ul}{\nwarrow}
\nc{\ur}{\nearrow}
\nc{\dl}{\swarrow}
\nc{\dr}{\searrow}
\nc{\dd}{\wedge}
\nc{\uu}{\vee}
\nc{\mm}{\ast}
\nc{\oo}{\circ}
\nc{\od}{\odot}
\nc{\oc}{\circledcirc}
\nc{\om}{\circledast}
\nc{\bk}{\blacksquare}
\nc{\rv}{\dashv}
\nc{\lv}{\vdash}

\nc{\vep}{\varepsilon}
\nc{\bin}[2]{ (_{\stackrel{\scs{#1}}{\scs{#2}}})}  %binomial coeff
\nc{\binc}[2]{(\!\! \begin{array}{c} \scs{#1}\\
    \scs{#2} \end{array}\!\!)}  %binomial coeff
\nc{\bincc}[2]{  ( {\scs{#1} \atop
    \vspace{-1cm}\scs{#2}} )}  %binomial coeff
\nc{\bs}{\bar{S}}
\nc{\la}{\longrightarrow}
\nc{\ot}{\otimes}
\nc{\rar}{\rightarrow}
\nc{\dar}{\downarrow}
\nc{\dap}[1]{\downarrow \rlap{$\scriptstyle{#1}$}}
\nc{\defeq}{\stackrel{\rm def}{=}}
\nc{\dis}[1]{\displaystyle{#1}}
\nc{\dotcup}{\ \displaystyle{\bigcup^\bullet}\ }
\nc{\hcm}{\ \hat{,}\ }
\nc{\hts}{\hat{\otimes}}
\nc{\hcirc}{\hat{\circ}}
\nc{\lleft}{[}
\nc{\lright}{]}
\nc{\curlyl}{\left \{ \begin{array}{c} {} \\ {} \end{array}
    \right .  \!\!\!\!\!\!\!}
\nc{\curlyr}{ \!\!\!\!\!\!\!
    \left . \begin{array}{c} {} \\ {} \end{array}
    \right \} }
\nc{\longmid}{\left | \begin{array}{c} {} \\ {} \end{array}
    \right . \!\!\!\!\!\!\!}
\nc{\ora}[1]{\stackrel{#1}{\rar}}
\nc{\ola}[1]{\stackrel{#1}{\la}}%${\Bbb Z}$
\nc{\scs}[1]{\scriptstyle{#1}} \nc{\mrm}[1]{{\rm #1}}
\nc{\dirlim}{\displaystyle{\lim_{\longrightarrow}}\,}
\nc{\invlim}{\displaystyle{\lim_{\longleftarrow}}\,}
\nc{\dislim}[1]{\displaystyle{\lim_{#1}}} \nc{\colim}{\mrm{colim}}
\nc{\mvp}{\vspace{0.3cm}} \nc{\tk}{^{(k)}} \nc{\tp}{^\prime}
\nc{\ttp}{^{\prime\prime}} \nc{\svp}{\vspace{2cm}}
\nc{\vp}{\vspace{8cm}}
\nc{\modg}[1]{\!<\!\!{#1}\!\!>}
\nc{\intg}[1]{F_C(#1)}
\nc{\lmodg}{\!<\!\!}
\nc{\rmodg}{\!\!>\!}
\nc{\cpi}{\widehat{\Pi}}
\nc{\sha}{{\mbox{\cyr X}}}  %used to be \cyr
\nc{\ssha}{{\mbox{\cyrs X}}} %sha as product
\nc{\tsha}{{\mbox{\cyrt X}}}
\nc{\shpr}{\diamond}    %Shuffle product
\nc{\labs}{\mid\!}
\nc{\rabs}{\!\mid}

\font\cyr=wncyr10
\font\cyrs=wncyr7
\font\cyrt=wncyr5

\nc{\ann}{\mrm{ann}}
\nc{\Aut}{\mrm{Aut}}
\nc{\can}{\mrm{can}}
\nc{\Cont}{\mrm{Cont}}
\nc{\rchar}{\mrm{char}}
\nc{\cok}{\mrm{coker}}
\nc{\dtf}{{R-{\rm tf}}}
\nc{\dtor}{{R-{\rm tor}}}

\nc{\Div}{{\mrm Div}}
\nc{\End}{\mrm{End}}
\nc{\Fil}{\mrm{Fil}}
\nc{\Fr}{\mrm{Fr}}
\nc{\Frob}{\mrm{Frob}}
\nc{\Gal}{\mrm{Gal}}
\nc{\GL}{\mrm{GL}}
\nc{\hsr}{\mrm{H}}
\nc{\hpol}{\mrm{HP}}
\nc{\im}{\mrm{im}}
\nc{\incl}{\mrm{incl}}
\nc{\length}{\mrm{length}}
\nc{\mchar}{\rm char}
\nc{\mpart}{\mrm{part}}
\nc{\ql}{{\QQ_\ell}}
\nc{\qp}{{\QQ_p}}
\nc{\rank}{\mrm{rank}}
\nc{\rcot}{\mrm{cot}}
\nc{\rdef}{\mrm{def}}
\nc{\rdiv}{{\rm div}}
\nc{\rtf}{{\rm tf}}
\nc{\rtor}{{\rm tor}}
\nc{\res}{\mrm{res}}
\nc{\SL}{\mrm{SL}}
\nc{\Spec}{\mrm{Spec}}
\nc{\tor}{\mrm{tor}}
\nc{\Tr}{\mrm{Tr}}
\nc{\tr}{\mrm{tr}}

\nc{\ab}{\mathbf{Ab}}
\nc{\DRB}{\mathbf{DRB}}\nc{\MDR}{\mathbf{MDR}}
\nc{\Alg}{{\mathbf{Alg}}}
\nc{\RBA}{{\mathbf{RBA}}}
\nc{\RB}{\mathbf{RB}_\lambda}
\nc{\bfk}{{\bf k}}
\nc{\bfone}{{\bf 1}}
\nc{\detail}{\marginpar{\bf More detail}
    \noindent{\bf Need more detail!}
    \svp}
\nc{\Diff}{\mathbf{Diff}}
\nc{\gap}{\marginpar{\bf Incomplete}\noindent{\bf Incomplete!!}
    \svp}
\nc{\FMod}{\mathbf{FMod}}
\nc{\Int}{\mathbf{Int}}
\nc{\Mon}{\mathbf{Mon}}
\nc{\Mult}{{\mathbf{Mlt}}}
\nc{\remarks}{\noindent{\bf Remarks: }}
\nc{\Rep}{\mathbf{Rep}}
\nc{\Rings}{\mathbf{Rings}}
\nc{\Sets}{\mathbf{Sets}}

\nc{\DD}{\mathbf{DD}}
\nc{\DT}{\mathbf{DT}}
\nc{\rt}{\triangleleft}
\nc{\lt}{\triangleright}
\nc{\Mod}{\mathbf{Mod}}
\nc{\Vect}{\mathbf{Vect}}
\nc{\VVect}{ \text{2Vect} }
\nc{\TTerm}{ \text{2Term} }
\nc{\BVVect}{\mathbf{2Vect}}
\nc{\BTTerm}{\mathbf{2Term}}
\nc{\LieAAlg}{\text{Lie2Alg}}\nc{\BLieAAlg}{\mathbf{Lie2Alg}}
\nc{\TTermL}{\text{2TermL}_\infty}\nc{\BTTermL}{\mathbf{2TermL}_\infty}
 \nc{\TTermRBL}{\text{2TermRB}L_\infty}
\nc{\BTTermRBL}{\mathbf{2TermRBL}_\infty}
\nc{\SLieAAlg}{\text{SLie2Alg}}\nc{\BSLieAAlg}{\mathbf{SLie2Alg}}
\nc{\XLie}{\text{XLie}}\nc{\BXLie}{\mathbf{XLie}}
\nc{\RBLie}{\text{RBLie}}\nc{\BRBLie}{\mathbf{RBLie}}
\nc{\RBLieAAlg}{\text{RBLie2Alg}}\nc{\BRBLieAAlg}{\mathbf{RBLie2Alg}}
\nc{\RBTTermL}{\text{RB2TermL_\infty}}\nc{\BRBTTermL}{\mathbf{RB2TermL_\infty}}
\nc{\SRBLieAAlg}{\text{SRBLie2Alg}}\nc{\BSRBLieAAlg}{\mathbf{SRBLie2Alg}}
\nc{\XRBLie}{\text{XRBLie}}\nc{\BXRBLie}{\mathbf{XRBLie}}
\nc{\SRBTTermL}{\mathbf{SRB2TermL_\infty}}
\nc{\BA}{{\mathbb A}}
\nc{\CC}{{\mathbb C}}
\nc{\EE}{{\mathbb E}}
\nc{\FF}{{\mathbb F}}
\nc{\GG}{{\mathbb G}}
\nc{\NN}{{\mathbb N}}
\nc{\PP}{{\mathbb P}}
\nc{\QQ}{{\mathbb Q}}
\nc{\RR}{{\mathbb R}}
\nc{\TT}{{\mathbb T}}
\nc{\VV}{{\mathbb V}}
\nc{\ZZ}{{\mathbb Z}}
\nc{\TP}{\widetilde{P}}
\nc{\lp}{\overline{P}}
\nc{\ld}{\overline{d}}
\nc{\SLT}{S}
\nc{\TTL}{T}

\nc{\cala}{{\mathcal A}}
\nc{\calc}{{\mathcal C}}
\nc{\cald}{\mathcal{D}}
\nc{\cale}{{\mathcal E}}
\nc{\calf}{{\mathcal F}}
\nc{\calg}{{\mathcal G}}
\nc{\calh}{{\mathcal H}}
\nc{\cali}{{\mathcal I}}
\nc{\call}{{\mathcal L}}
\nc{\calm}{{\mathcal M}}
\nc{\caln}{{\mathcal N}}
\nc{\calo}{{\mathcal O}}
\nc{\calp}{{\mathcal P}}
\nc{\calr}{{\mathcal R}}
\nc{\cals}{{\mathcal S}}
\nc{\calt}{{\mathcal T}}
\nc{\calu}{{\mathcal U}}
\nc{\calv}{{\mathcal V}}
\nc{\calw}{{\mathcal W}}
\nc{\calx}{{\mathcal X}}
\nc{\CA}{\mathcal{A}}

\nc{\fraka}{{\mathfrak a}}
\nc{\frakb}{\mathfrak{b}}
\nc{\frakB}{{\frak B}} \nc{\frakm}{{\frak
m}} \nc{\frakM}{{\frak M}}
\nc{\frakp}{{\frak p}}
\nc{\frakS}{{\frak S}}
\nc{\frakA}{{\frak A}} \nc{\frakx}{{\frakx}}

\nc{\Dif}{\mathbf{Dif}}
\nc{\Difl}{\mathbf{Dif}_\lambda}
\nc{\Lie}{\mathbf{Lie}}
\nc{\DifLie}{\mathbf{DifLie}}
\nc{\XAlg}{\mathbf{XAlg}}
\nc{\XDif}{\mathbf{XDif}}
\nc{\XRB}{\mathbf{XRB}}
\nc{\XDRB}{\mathbf{XDRB}}
\nc{\XDD}{\mathbf{XDD}}
\nc{\XDT}{\mathbf{XDT}}
\nc{\XP}{\mathbf{XP}}
\nc{\SP}{\mathbf{2-P}}
\nc{\SSAlg}{\mathbf{2-Alg}}
\nc{\SAlg}{\mathbf{2\text{-}SAlg}}
\nc{\SDif}{\mathbf{2-SDif}}
\nc{\SRB}{\mathbf{2-SRB}}
\nc{\SDRB}{\mathbf{2-SDRB}}
\nc{\SDD}{\mathbf{2-SDD}}
\nc{\SDT}{\mathbf{2-SDT}}
\nc{\SLie}{\mathbf{2-SLie}}
\nc{\SDifLie}{\mathbf{2-SDifLie}}
\nc{\SF}{\mathbf{S\calf}}
\nc{\XF}{\mathbf{X\calf}}
\nc{\C}{\mathbf{C}}
\nc{\A}{\mathbf{A}}
\nc{\B}{\mathbf{B}}
\nc{\T}{\mathbf{T}}
\nc{\SDP}{\mathbf{2-SDP}}
\nc{\SRBP}{\mathbf{2-SRBP}}
\nc{\XDP}{\mathbf{XDP}}
\nc{\XRBP}{\mathbf{XRBP}}
\newtheorem{thm}{Theorem}[section]

\newtheorem{cor}[thm]{Corollary}
\newtheorem{pro}[thm]{Proposition}
\newtheorem{ex}[thm]{Example}

\newtheorem{defi}[thm]{Definition}

\newcommand {\emptycomment}[1]{} %to remove paragraphs

\newcommand{\lon }{\,\rightarrow\,}
\newcommand{\be }{\begin{equation}}
\newcommand{\ee }{\end{equation}}

\newcommand{\huaR}{\mathcal{R}}

\newcommand{\huaG}{\mathcal{G}}

\newcommand{\huaD}{\mathcal{D}}

\newcommand{\huaO}{\mathcal{O}}

\newcommand{\CE}{\mathsf{CE}}

\newcommand{\frkC}{\mathfrak C}

\newcommand{\frkR}{\mathfrak R}

\newcommand{\br}[1]{   [ \cdot,    \cdot  ]   }

\newcommand{\g}{\mathfrak g}
\newcommand{\h}{\mathfrak h}

\newcommand{\dM}{\mathrm{d}}

\newcommand{\Hom}{\mathrm{Hom}}
\newcommand{\Der}{\mathrm{Der}}

\newcommand{\gl}{\mathfrak {gl}}

\newcommand{\Ker}{\mathrm{ker}}
\newcommand{\ad}{\mathrm{ad}}

\begin{document}
\title[Rota-Baxter Lie $2$-algebras]{Rota-Baxter Lie $2$-algebras}

\author{Shilong Zhang}
\address{College of Science, Northwest A\&F University, Yangling 712100, Shaanxi, China }
\email{shlzhang11@163.com}

\author{Jiefeng Liu$^{\ast}$}
\address{School of Mathematics and Statistics, Northeast Normal University,\\
 Changchun 130024, Jilin, China }
\email{liujf534@nenu.edu.cn}
\thanks{$^{\ast}$ the corresponding author}

\vspace{-5mm}

\begin{abstract}
In this paper, we introduce the notion of Rota-Baxter Lie $2$-algebras, which is a categorification of Rota-Baxter Lie algebras. We prove that the category of Rota-Baxter Lie
$2$-algebras and the category of $2$-term Rota-Baxter $L_\infty$-algebras are equivalent.
We introduce the notion of  a crossed module of Rota-Baxter Lie algebras and show that there is a one-to-one correspondence between strict $2$-term Rota-Baxter $L_\infty$-algebras and crossed modules of Rota-Baxter Lie algebras.  We give the construction of crossed modules of  Lie algebras from crossed modules of Rota-Baxter Lie algebras.
\end{abstract}

\keywords{Rota-Baxter Lie $2$-algebra, $2$-term Rota-Baxter $L_\infty$-algebra, crossed module}
\footnotetext{{\it{MSC}}: 17B38, 18N25}

\maketitle

\tableofcontents

\setcounter{section}{0}

\allowdisplaybreaks
\section{Introduction}
The concept of Rota-Baxter operators on associative algebras was introduced in the 1960s by G. Baxter \cite{Ba} and G.-C. Rota \cite{Rota} in their study of
fluctuation theory in probability and combinatorics. Recently it has found many applications, including in Connes-Kreimer's algebraic approach to the renormalization in perturbative quantum field theory \cite{CK}. In the Lie algebra context, a Rota-Baxter operator was introduced independently in the 1980s as the operator form of the classical Yang-Baxter equation, named after the physicists C.-N. Yang and R. Baxter \cite{Bax,Yang}, whereas the classical Yang-Baxter equation plays important roles in many fields in mathematics and mathematical physics such as integrable systems and quantum groups \cite{CP,Semonov-Tian-Shansky}. A Lie algebra equipped with a Rota-Baxter operator is called a Rota-Baxter Lie algebra.  Recently, cohomologies, deformations and extensions of Rota-Baxter Lie algebras are studied in \cite{JS,LST,TBGS}. See \cite{Gub} fore details on Rota-Baxter operators.

Motivated by the study of string theory, people pay more attention to higher categorical structures. One way to obtain higher categorical structures is by categorifying existing mathematical concepts. One of the simplest higher structures is a $2$-vector space, which is a categorification of a vector space. If we further put Lie algebra structures on $2$-vector spaces, then we obtain Lie $2$-algebras \cite{BC}.  $L_\infty$-algebras, sometimes called strongly homotopy Lie algebras, were introduced \cite{LS} as a model for  Lie algebras that satisfy Jacobi identity up to all higher homotopies. It is well-known that the category of Lie $2$-algebras is equivalent to the category of $2$-term $L_\infty$-algebras. The structure of a $2$-term $L_\infty$-algebra appears in many areas such as string theory \cite{BR}, higher symplectic geometry \cite{BHC}, and Courant algebroids \cite{Roy}.

In this paper, we provide the categorification of Rota-Baxter Lie algebras, which we call Rota-Baxter Lie $2$-algebra. Rota-Baxter operators (more generally, $\huaO$-operators) on $2$-term $L_\infty$-algebras were first introduced in \cite{Sh} as a tool to study $2$-graded classical Yang-Baxter equations, which could naturally generate examples of Lie 2-bialgebras \cite{BSZ}. Soon afterwards, Rota-Baxter operators on $L_\infty$-algebras were given and studied   in \cite{LST}. We prove that the category of Rota-Baxter Lie 2-algebras and category of  $2$-term Rota-Baxter $L_\infty$-algebras are equivalent. Here  a $2$-term Rota-Baxter $L_\infty$-algebra consists of a  $2$-term $L_\infty$-algebra and a Rota-Baxter operator on it. The notion of crossed modules of Rota-Baxter Lie algebras is also introduced and we prove that there is a one-to-one correspondence between strict $2$-term Rota-Baxter $L_\infty$-algebras and crossed modules of Rota-Baxter Lie algebras. We show that a crossed module of Rota-Baxter Lie algebras gives a crossed module of pre-Lie algebras and thus gives a crossed module of Lie algebras naturally.

The paper is organized as follows. In Section \ref{sec:Prel}, we recall Rota-Baxter Lie algebras and their representations, $2$-vector spaces and $2$-term chain complexes. In Section \ref{sec:Lie2Alg}, we first give the notion of Rota-Baxter Lie $2$-algebras, which is the categorification of Rota-Baxter Lie algebras. Then we introduce the category of Rota-Baxter Lie 2-algebras and the category of $2$-term Rota-Baxter $L_\infty$-algebras and show that they are equivalent. In Section \ref{sec:sRBLie},  we introduce the notion of crossed modules of Rota-Baxter Lie algebras  and show that there is a one-to-one correspondence between strict $2$-term Rota-Baxter $L_\infty$-algebras and crossed modules of Rota-Baxter Lie algebras. We show that the underlying algebraic structure of  a crossed module of Rota-Baxter Lie algebras is a crossed module of pre-Lie algebras and then a new crossed module of Lie algebras is constructed.

In this paper, all the vector spaces are over algebraically closed field $\mathbb K$ of characteristic $0$, and finite dimensional.

\section{Preliminaries}\label{sec:Prel}
\subsection{Rota-Baxter Lie algebras and their representations}
\begin{defi}
Let $(\g,[\cdot,\cdot]_\g)$ be a Lie algebra.
A linear operator $R:\g\longrightarrow \g$ is called a {\bf Rota-Baxter operator } if
\begin{equation*}
 [R(x),R(y)]_\g=R\big([R(x),y]_\g+ [x,R(y)]_\g \big), \quad \forall x, y \in \g.
\end{equation*}
Moreover, a Lie algebra $(\g,[\cdot,\cdot]_\g)$ with a Rota-Baxter operator $R$ is
called a {\bf Rota-Baxter Lie algebra}. We denote it by $(\g,[\cdot,\cdot]_\g,R)$.
\end{defi}

%As a generalization of a Rota-Baxter operator on a Lie algebra to an arbitrary representation, Kupershmidt introduce the notion of an $\mathcal{O}$-operator $($ also called a Kupershmidt operator$)$\cite{Ku}:

%\begin{defi}
%Let $(\g,[\cdot,\cdot]_\g)$ be a Lie algebra and $(V;\rho)$ be a representation of $(\g,[\cdot,\cdot]_\g)$. A linear map $T:V\rar \g$ is called an {\bf $\mathcal{O}$-operator } on $\g$ associated to the representation $(V;\rho)$ if $T$ satisfies
%\begin{equation*}
% [T(u),T(v)]_\g=T\big(\rho(T(u))v-\rho(T(v))u \big), \quad \forall u, v \in V.
%\end{equation*}
%\end{defi}

%Then a Rota-Baxter operator $R$ on a Lie algebra $\g$ is an $\mathcal{O}$-operator with respect to the adjoint representation $(\g;\text{ad})$:
%\begin{equation*}
% [R(x),R(y)]_\g=R\big(\text{ad}_{R(x)}y-\text{ad}_{R(y)}x \big), \quad \forall x, y \in \g.
%\end{equation*}

\begin{defi}  A {\bf pre-Lie algebra} is a pair $(\g,\ast_\g)$, where $\g$ is a vector space and  $\ast_\g:\g\otimes \g\longrightarrow \g$ is a bilinear multiplication satisfying that for all $x,y,z\in \g$, the associator
$(x,y,z)=(x\ast_\g y)\ast_\g z-x\ast_\g(y\ast_\g z)$ is symmetric in $x,y$,
i.e.
$$(x,y,z)=(y,x,z),\;\;{\rm or}\;\;{\rm
equivalently,}\;\;(x\ast_\g y)\ast_\g z-x\ast_\g(y\ast_\g z)=(y\ast_\g x)\ast_\g z-y\ast_\g(x\ast_\g z).$$
\end{defi}

Let $(\g,\ast_\g)$ be a pre-Lie algebra. The commutator
$ [x,y]_\g=x\ast_\g y-y\ast_\g x$ defines a Lie algebra structure on $\g$,
which is called the {\bf sub-adjacent Lie algebra} of $(\g,\ast_\g)$ and denoted by $\g^c$. Furthermore, $L:\g\rightarrow
\gl(\g)$ defined by
\begin{equation}\label{eq:defiLpreLie}
L_xy=x\ast_\g y,\quad \forall x,y\in \g
\end{equation}
 gives a representation of $\g^c$ on $\g$. See \cite{Pre-lie algebra in
geometry} for more details.

The following proposition is well-known.
\begin{pro}
 Let $(\g,[\cdot,\cdot]_\g)$ be a Lie algebra and $R:\g\longrightarrow
\g$ a Rota-Baxter operator. Define a new operation on $\g$ by
$$x\ast y=[R(x),y]_\g.$$
Then $(\g,\ast)$ is a pre-Lie algebra and $R$ is a homomorphism from the sub-adjacent Lie algebras $(\g,[-,-]_R)$ to $(\g,[-,-]_\g)$, where $[x,y]_R=x\ast y-y\ast x$.
\end{pro}

\begin{defi}
  Let $(\g,[\cdot,\cdot]_\g,R)$ and $(\h,[\cdot,\cdot]_\h,S)$ be two Rota-Baxter Lie algebras. A {\bf homomorphism} is a linear map  $\phi:\g\rightarrow\h$ such that $\phi$ a Lie algebra homomorphism and satisfies $\phi\circ R=S\circ \phi.$

\end{defi}

A {\bf Rota-Baxter Lie subalgebra} (resp., {\bf Rota-Baxter Lie ideal}) of a Rota-Baxter Lie algebra $(\g, R)$ is a Lie subalgebra (resp., a Lie ideal) $I$ of $\g$ such that $R(I)\subseteq I$.

Let $f:(\g,R)\rar (\h,S)$ be a Rota-Baxter Lie algebra homomorphism. Then  $\Ker~f$ is a Rota-Baxter Lie ideal of the Rota-Baxter Lie algebra $(\g,R)$.

\emptycomment{A {\bf derivation} on a Rota-Baxter Lie algebra $(\g,[\cdot,\cdot]_{\g},R)$ is a linear map  $ D \in\Hom(\g,\g)$   such that $D$ is a derivation on  the Lie algebra $\g$ and $D\circ R=R\circ D$. Denote the set of derivations on Rota-Baxter Lie algebra $(\g,[\cdot,\cdot]_{\g},R)$ by $\Der_{\rm RB}(\g)$.}

\begin{defi}{\rm (\cite{JS} )}
A {\bf representation of a Rota-Baxter Lie algebra} $(\g,[\cdot,\cdot]_\g, R)$ on a vector space $V$ with respect to a linear map $\huaR\in\gl(V)$ is a representation $\rho$ of the Lie algebra $\g$ on $V$, satisfying
\begin{equation}\label{eq:rep-RB}
\rho(R(x))\circ \huaR =\huaR\circ\rho(R(x)) +\huaR\circ \rho(x)\circ \huaR, \quad \forall x\in\g.
\end{equation}

\end{defi}

Denote a representation by $(V; \rho,\huaR)$.

\begin{ex}{\rm
Let $(\g,[\cdot,\cdot]_{\g},R)$  be a  Rota-Baxter Lie algebra.
Then $(\g;\ad,R)$ is a representation, which is called the {\bf adjoint representation} of $(\g,[\cdot,\cdot]_{\g}, R)$.
}
\end{ex}

\begin{pro}
 Let $(V; \huaR, \rho)$ be a representation of a  Rota-Baxter Lie algebra $(\g,[\cdot,\cdot]_{\g},R)$. Then $(V^*; \rho^*,-\huaR^*)$ is   also a representation of   $(\g,[\cdot,\cdot]_{\g},R)$, which is called the {\bf dual representation}.
\end{pro}
\begin{ex}{\rm
Let $(\g,[\cdot,\cdot]_{\g},R)$  be a  Rota-Baxter Lie algebra. Then $(\g^*;\ad^*,-R^*)$ is a representation of $(\g,[\cdot,\cdot]_{\g},R)$, which is called the {\bf coadjoint representation}.
}
\end{ex}

\begin{pro}
Let $(V;\rho, \huaR)$ be a representation of a Rota-Baxter Lie algebra $(\g,[\cdot,\cdot]_{\g}, R)$. Then $(\g\oplus V, [\cdot,\cdot]_{\ltimes}, \frak{R})$ is a Rota-Baxter Lie algebra, where $[\cdot,\cdot]_{\ltimes}$ is the semidirect product Lie bracket given by
\begin{equation*}
[x+u,y+v]_{\ltimes}=[x,y]_\g+\rho(x)v-\rho(y)u,\quad \forall x,y\in\g, u,v\in V,
\end{equation*}
and $\frak{R}: \g\oplus V\rightarrow\g\oplus V$ is a linear map given by
\begin{equation*}
\frak{R}(x+u)=R(x)+\huaR(u),\quad \forall x\in\g, u\in V.
\end{equation*}
\end{pro}

\subsection{$2$-vector spaces}
Let $\Vect$ be the category of vector spaces. Vector spaces can be categorified to $2$-vector spaces. A good introduction for this subject is \cite{BC}.
\begin{defn}
A {\bf $2$-vector space} is an internal category in the category Vect.
\end{defn}

Thus, a $2$-vector space $V=(V_1, V_0, s, t, i,\circ)$ is a category with a vector space of objects $V_0$ and a vector space of morphisms $V_1$, such that the source and target maps $s, t:V_1\rightarrow V_0$, the identity-assigning map $i:V_0\rightarrow V_1$, and the composition map $\circ:V_1\times_{V_0}V_1\rightarrow V_1$ are linear.

Given a morphism $f:x\rightarrow y\in V_1$, define the {\bf arrow part} of $f$, denoted as $\vec{f}$, by
$$\vec{f}= f-i(x).$$
Furthermore, we identify $f:x\rightarrow y$ with the ordered pair $(x,\vec{f})$. It was shown in \cite{BC} that the
composition map $\circ: V_1\times V_1 \rar V_1$ is uniquely determined by
\begin{equation}
f\circ g=(x,\vec{f}+\vec{g}),\quad f=(x,\vec{f}),g=(y,\vec{g})\in V_1.\label{eq:ppcirc}
\end{equation}
Thus the structure of a $2$-vector space is completely determined by the vector spaces $V_0$ and $V_1$ together with the source, target and identity-assigning maps.

Let $V$ and $W$ be two $2$-vector spaces. Recall that a {\bf linear functor} $F:V\rightarrow W$ is an internal functor in Vect.

Let $\VVect$ denote the category consisting of $2$-vector spaces and linear functors between them. There is a category, denoted as $\TTerm$, whose objects are $2$-term chain complexes and whose morphisms are chain maps.

It is well known that the categories $\VVect$ and $\TTerm$ are equivalent. Roughly speaking, given a $2$-vector space $V=(V_1, V_0, s, t, i,\circ)$,
$$\Ker(s)\stackrel{t}{\rightarrow} V_0$$
is a $2$-term complex. Conversely, the $2$-term complex of vectors $C_1\stackrel{d}{\rightarrow} C_0$ gives a $2$-vector space of which the set of objects is $C_0$, the set of morphism is $C_0\oplus C_1$, the identity-assigning map is given by  $i(x)=(x,0)$ for $x\in C_0$, the source map $s$ is given by $s(x,\vec{f})=x$ and the target map $t$ is given by $t(x,\vec{f})=x+d\vec{f}$ for $(x,\vec{f})\in C_0\oplus C_1$.

\begin{defn}
\begin{enumerate}
\item
Given two linear functors $F, G:V\rar W$ between $2$-vector spaces, a {\bf linear natural transformation} $\alpha:F\Rightarrow G$ is a natural transformation in Vect.
\item
Given two chain maps $\varphi, \psi:C\rar C'$ of $2$-term chain complexes, a {\bf chain homotopy} $\tau:\varphi\Rightarrow \psi$ is a map $\tau:C_0\rar C_1'$ satisfying $d'\tau=\psi_0-\varphi_0$ and $\tau d=\psi_1-\varphi_1$.
\end{enumerate}
\end{defn}

Let $\BVVect$ denote the $2$-category of $2$-vector spaces, linear functors and linear natural transformations. Also let $\BTTerm$ be the $2$-category of $2$-term chain complexes, chain maps, and chain homotopies.

Furthermore, we have
\begin{prop}(\cite{BC})
The $2$-category 2Vect is $2$-equivalent to the $2$-category 2Term.
\label{prop:2BVect2Term}
\end{prop}

\section{Rota-Baxter Lie $2$-algebras and $2$-term Rota-Baxter $L_\infty$-algebras}\label{sec:Lie2Alg}

In this section, we first introduce the notion of a Rota-Baxter Lie $2$-algebra which is a  Lie $2$-algebra with a linear functor which satisfies the Rota-Baxter identity up to a natural isomorphism. Then we introduce the notion of a $2$-term Rota-Baxter $L_\infty$-algebra. Finally, we show that the category of Rota-Baxter Lie
$2$-algebras and the category of $2$-term Rota-Baxter $L_\infty$-algebras are equivalent.
\subsection{Rota-Baxter Lie $2$-algebras }
\label{subsec:Lie2Alg}

We begin by reviewing the concept of a Lie $2$-algebra given in \cite{BC}.
\begin{defi}
\begin{enumerate}
\item
A {\bf Lie $2$-algebra} is a $2$-vector space $ L$ together with a skew-symmetric bilinear functor $[\cdot,\cdot]:L\times L\rar L$ and a completely antisymmetric trilinear natural isomorphism, the {\bf Jacobiator},
$$J_{x, y, z}:[[x,y],z]\rar[x,[y,z]]+[[x,z],y],$$
satisfying the identity:
\begin{eqnarray*}
&&([w,J_{x,y,z}]+1)([J_{w,y,z},x]+1)
(J_{[w,y],x,z}+J_{w,[x,y],z})[J_{w,x,y},z]
\\
&=&(J_{w,[x,z],y}+J_{[w,z],x,y}+J_{w,x,[y,z]})
([J_{w,x,z},y]+1)J_{[w,x],y,z}.
\end{eqnarray*}
A Lie $2$-algebra is called {\bf strict} if the Jacobiator is the identity isomorphism.
\item
Given two Lie $2$-algebras $ L$ and $ L'$, a {\bf homomorphism} $F=(F_0,F_1,F_2): L\rar L'$ consists of a linear functor $(F_0,F_1)$ from the underlying $2$-vector space of $ L$ to that of $L'$, and a skew-symmetric bilinear natural transformation
$$F_2[x,y]:[F_0(x),F_0(y)]\rar F_0[x,y]$$
satisfying
$$(F_1(J_{x,y,z}))F_2[F_2,1]
=(F_2+F_2)([1,F_2]+[F_2,1])J_{F_0(x),F_0(y),F_0(z)}.$$
\end{enumerate}
\label{defn:ssLie}
\end{defi}

In the following, we give the main definition in this paper.
\begin{defn}
A {\bf  Rota-Baxter Lie $2$-algebra} is a triple $((L,[\cdot,\cdot]),P,\huaR)$, where $(L,[\cdot,\cdot])$ is a Lie $2$-algebra, $P=(P_0,P_1):L \rar L$ is a linear functor  and for $x,y\in L$, $\huaR_{x,y}$ is an antisymmetric bilinear natural isomorphism given by
$$\huaR_{x, y}:[P_{0}(x),P_{0}(y)]\rar P_{0}[P_{0}(x),y]+P_{0}[x,P_{0}(y)],$$
such that the following Rota-Baxter relation is satisfied,
{\small{
\begin{eqnarray}
&&\left(1+P_{1}[\huaR_{x,y},1_{z}]\right)\left(1+P_{1}J_{P_{0}(x),z,P_{0}(y)}\right)\left(\huaR_{x,[P_{0}(y),z]}+\huaR_{x,[y,P_{0}(z)]}
+\huaR_{[x,P_{0}(z)],y}+\huaR_{[P_{0}(x),z],y}\right)\notag\\
\label{eq:RBcoh}&&\left([1_{P_{0}(x)},\huaR_{y,z}]+[\huaR_{x,z},1_{P_{0}(y)}]\right)J_{P_{0}(x),P_{0}(y),P_{0}(z)}\\
\notag&=&\left(1+P_{1}[\huaR_{x,z},1_{y}]+P_{1}[1_{x},\huaR_{y,z}]\right)\left(1+P_{1}J_{P_{0}(x),y,P_{0}(z)} +P_{1}J_{x,P_{0}(y),P_{0}(z)}\right)\\
&&\left(\huaR_{[P_{0}(x),y], z} +\huaR_{[x,P_{0}(y)], z}\right)[\huaR_{x,y},1_{P_{0}(z)}],
\notag
\end{eqnarray}
}}
which can be showed as the following commutative diagram

{\tiny{\begin{equation*}
\begin{array}{c}
\xymatrix{
& { [[P_{0}(x),P_{0}(y)],P_{0}(z)]} \ar[dl]_{J_{P_{0}(x),P_{0}(y),P_{0}(z)}} \ar[rd]^{1}&\\
[P_{0}(x),[ P_{0}(y) ,P_{0}(z)]]+[[P_{0}(x),P_{0}(z)],P_{0}(y)]
\ar[d]_{[1_{P_{0}(x)},\huaR_{y,z}]+[\huaR_{x,z},1_{P_{0}(y)}]}
&&[[P_{0}(x),P_{0}(y)],P_{0}(z)]\ar[d]^{[\huaR_{x,y},1_{P_{0}(z)}]}\\
A
\ar[d]_{\huaR_{x,[P_{0}(y),z]}+\huaR_{x,[y,P_{0}(z)]}
+\huaR_{[x,P_{0}(z)],y}+\huaR_{[P_{0}(x),z],y}}
&& [P_{0}[P_{0}(x),y],P_{0}(z)]+ [P_{0}[x,P_{0}(y)],P_{0}(z)] \ar[d]^{ \huaR_{[P_{0}(x),y], z} +\huaR_{[x,P_{0}(y)], z}}\\
{B
}\ar[d]_{1+P_{1}J_{P_{0}(x),z,P_{0}(y)}}
&&{D}
\ar[d]^{1+P_{1}J_{P_{0}(x),y,P_{0}(z)}+P_{1}J_{x,P_{0}(y),P_{0}(z)} }\\
{C}\ar[dr]_{1+P_{1}[\huaR_{x,y},1_{z}]}&&{E}
\ar[dl]^{\quad 1+P_{1}[\huaR_{x,z},1_{y}]+P_{1}[1_{x},\huaR_{y,z}]}
\\
&F
 &
}
\end{array}
\end{equation*}
}}
where \begin{eqnarray*}
 A&=&[P_{0}(x),P_{0}[P_{0}(y),z]]+[P_{0}(x),P_{0}[y,P_{0}(z)]]+[P_{0}[x,P_{0}(z)],P_{0}(y)]+[P_{0}[P_{0}(x),z],P_{0}(y)];\\
B&=&P_{0}[P_{0}(x),[ P_{0}(y) ,z]] +P_{0}[x,P_{0}[ P_{0}(y) ,z]]
+P_{0}[P_{0}(x), [ y,P_{0}(z)]]+P_{0}[x,P_{0}[ y,P_{0}(z)]]\\
&&+P_{0}[P_{0}[x,P_{0}(z)],y]+P_{0}[[x,P_{0}(z)],P_{0}(y)]
 + P_{0}[P_{0}[P_{0}(x),z],y]+P_{0}[[P_{0}(x),z],P_{0}(y)];\\
 C&=&P_{0}[P_{0}(x),[ P_{0}(y) ,z]] +P_{0}[x,P_{0}[ P_{0}(y) ,z]]+P_{0}[P_{0}(x), [ y,P_{0}(z)]]+P_{0}[x,P_{0}[ y,P_{0}(z)]]\\
&&+P_{0}[P_{0}[x,P_{0}(z)],y]+P_{0}[[x,P_{0}(z)],P_{0}(y)]+ P_{0}[P_{0}[P_{0}(x),z],y]+P_{0}[P_{0}(x),[z,P_{0}(y)]]\\
&&+P_{0}[[P_{0}(x),P_{0}(y)],z];\\
D&=&P_{0}[P_{0}[P_{0}(x),y],z]+P_{0}[P_{0}[x,P_{0}(y)],z]+P_{0}[[P_{0}(x),y],P_{0}(z)]+ P_{0}[[x,P_{0}(y)],P_{0}(z)];\\
E&=&P_{0}[P_{0}[P_{0}(x),y],z]+P_{0}[P_{0}[x,P_{0}(y)],z]+ P_{0}[P_{0}(x),[ y ,P_{0}(z)]]+P_{0}[[x,P_{0}(z)],P_{0}(y)]\\
&&+P_{0}[[P_{0}(x),P_{0}(z)],y]+ P_{0}[x,[P_{0}(y),P_{0}(z)]];\\
F&=&P_{0}[P_{0}[P_{0}(x),y],z]+P_{0}[P_{0}[x,P_{0}(y)],z] +P_{0}[P_{0}(x),[ y ,P_{0}(z)]]+P_{0}[[x,P_{0}(z)],P_{0}(y)]\\
&&+P_{0}[P_{0}[P_{0}(x),z],y]+P_{0}[P_{0}[x,P_{0}(z)],y]+P_{0}[x,P_{0}[P_{0}(y),z]]+P_{0}[x,P_{0}[y,P_{0}(z)]].
\end{eqnarray*}
A Rota-Baxter Lie $2$-algebra is called {\bf strict} if $(L,[\cdot,\cdot])$ is a strict Lie $2$-algebra and the natural isomorphism $R$ is the identity isomorphism.
\end{defn}

\begin{defn}\label{defn:RBLh}
Let $((L,[\cdot,\cdot]),P,\huaR)$ and $((L',[\cdot,\cdot]),P',\huaR')$ be two Rota-Baxter Lie $2$-algebras. A {\bf homomorphism of Rota-Baxter Lie $2$-algebras} $F:L\rar L'$ consists of a homomorphism of  Lie $2$-algebras $(F_0,F_1,F_2):L\rar L'$ and a natural linear transformation
$$F_3(x):P_0'(F_0(x))\rar F_0(P_0(x))$$
such the following equation holds
\begin{eqnarray}
&&(F_3[P_0(x),y]+F_3[x,P_0(y)])(P_1'F_2(P_0(x),y) +P_1'F_2(x,P_0(y)))\notag\\
&&(P_1'[F_3(x),1_{F_0(y)}]  +P_1'[1_{F_0(x)},F_3(y)])\huaR_{F_0(x),F_0(y)}\notag\\
&=&F_1(\huaR_{x,y})F_2(P_0(x),P_0(y))[F_3(x),F_3(y)],\label{eq:RBcohm}
\end{eqnarray}
or, in terms of diagram,
{\small\begin{equation*}
\begin{array}{c}
\xymatrix{
[P_0'(F_0(x)),P_0'(F_0(y))]\ar[d]_{[F_3(x),F_3(y)]} \ar[r]^{\quad\quad \huaR_{F_0(x),F_0(y)}\quad}
& P_0'[P_0'(F_0(x)),F_0(y)]\atop +P_0'[F_0(x),P_0'(F_0(y))]\ar[r]^{P_1'[F_3(x),1_{F_0(y)}]\atop +P_1'[1_{F_0(x)},F_3(y)]}
&P_0'[F_0(P_0(x)),F_0(y)]\atop +P_0'[F_0(x),F_0(P_0(y))]\ar[d]^{P_1'F_2(P_0(x),y) +P_1'F_2(x,P_0(y))}\\
[F_0(P_0(x)),F_0(P_0(y))]\ar[d]_{F_2(P_0(x),P_0(y))}&&P_0'(F_0[P_0(x),y]) +P_0'(F_0[x, P_0(y)]) \ar[d]^{F_3[P_0(x),y]+F_3[x,P_0(y)]} \\
F_0[P_0(x),P_0(y)]\ar[rr]_{F_1(\huaR_{x,y})}&&F_0P_0[P_0(x),y]+F_0P_0[x,P_0(y)]
.}
\end{array}
\end{equation*}}
\end{defn}
\vspace{3mm}

Let $(L,P,\huaR)$, $(L',P',\huaR')$ and $(L'',P'',\huaR'')$ be Rota-Baxter Lie $2$-algebras. Let $F:(L,P,\huaR)\rar (L',P',\huaR')$ and $G:(L',P',\huaR')\rar (L'',P'',\huaR'')$ be homomorphisms of Rota-Baxter Lie $2$-algebras. We define the composite functor $G\circ F:(L,P,\huaR)\rar (L'',P'',\huaR'')$ to be the usual composite of the underlying $2$-vector space functor:
$L\stackrel{F}{\longrightarrow}L'\stackrel{G}{\longrightarrow}L'',$
while letting $(G\circ F)_2$ and  $(G\circ F)_3$ be defined as the following composite
\begin{eqnarray*}
(G\circ F)_2[(G\circ F)_0(x),(G\circ F)_0(y)]&=&(G\circ F_2)\left(G_2[G_0(F_0(x)), G_0(F_0(y))]\right)=(G\circ F)_0[x,y],\\
(G\circ F)_3\left(P''((G\circ F)_0(x))\right)&=&(G\circ F_3)\left(G_3(P''(G_0(F_0(x))))\right)=(G\circ F)_0\left(P(x)\right),
\end{eqnarray*}
where $G\circ F_2$( resp. $G\circ F_3$) is the result of whiskering the functor $G$ by the natural transformation $F_2$( resp. $F_3$). The identity homomorphism $1_{L}$ has the identity functor as its underlying functor, together identity natural transformations $(1_{L})_2$ and $(1_{L})_3$. It is straightforward to obtain
\begin{prop}
  There is a category, which we denote by {\bf RBLie2Alg}, with Rota-Baxter Lie $2$-algebras as objects, Rota-Baxter Lie $2$-algebra homomorphisms as morphisms.
\end{prop}

%\tableofcontents
\subsection{A $2$-term Rota-Baxter $L_\infty$-algebras}
The notion of an $L_\infty$-algebra was introduced by Stasheff in \cite{LS}.  We begin by reviewing the concept of a $2$-term $L_\infty$-algebra.
\begin{defi}\label{defi:2-term}
  A {\bf $2$-term $L_\infty$-algebra} on a graded vector space $\huaG=\g_0\oplus \g_1$ consists of the following data:
\begin{itemize}
\item[$\bullet$] a complex of vector spaces: $\g_{1}\stackrel{l_1}{\longrightarrow}\g_0,$

\item[$\bullet$] a skew-symmetric bilinear map $l_2:\g_{i}\otimes \g_{j}\longrightarrow
\g_{i+j}$, where  $0\leq i+j\leq1$, which we denote more suggestively as $[\cdot,\cdot]$,

\item[$\bullet$] a  skew-symmetric trilinear map $l_3:\wedge^3\g_0\longrightarrow
\g_{1}$,
   \end{itemize}
   such that for all $x_i,x,y,z\in \g_0$ and $u,v\in \g_{1}$, the following equalities are satisfied:
\begin{itemize}
\item[$\rm(a)$] $l_1 l_2(x,u)=l_2(x,l_1(u)),\quad l_2(l_1( u),v)=l_2( u,l_1(v)),$
\item[$\rm(b)$]$l_1 l_3(x,y,z)=l_2(x,l_2(y,z))+l_2(z,l_2(x,y))+l_2(y,l_2(z,x)),$
\item[$\rm(c)$]$ l_3(x,y,l_1(u))=l_2(x,l_2(y,u))+l_2(u,l_2(x,y))+l_2(y,l_2(u,x)),$
\item[$\rm(d)$] the Jacobiator identity:
\begin{eqnarray*}
&&\sum_{i=1}^4(-1)^{i+1}l_2(x_i,l_3(x_1,\cdots,\hat{x_i}),\cdots,x_4)\\
&&+\sum_{i<j}l_3(l_2(x_i,x_j),x_1,\cdots,\hat{x_i},\cdots,\hat{x_j},\cdots,x_4)=0.\end{eqnarray*}
   \end{itemize}
\label{defi:2TL}
\end{defi}
We usually denote a $2$-term $L_\infty$-algebra by $(\g_{1},\g_0,l_1,l_2,l_3)$, or simply by $\huaG$. A $2$-term $L_\infty$-algebra is called  {\bf strict} if $l_3=0$.
\begin{defi}
Let $\huaG=(\g_{1},\g_0,l_1,l_2,l_3)$ and $\huaG'=(\g_{1}',\g_0',l_1',l_2',l_3')$ be $2$-term $L_\infty$-algebras. An {\bf $L_\infty$-homomorphism} $\phi:\huaG\rar\huaG'$ consists of:
\begin{itemize}
  \item a chain map $\phi :\huaG \to \huaG'$ which consists of
  linear maps $\phi_0 : \g_0 \to \g'_0$ and
              $\phi_1 : \g_1 \to \g'_1$ preserving the differential: $l_1'\phi_1=\phi_0l_1$,
  \item a skew-symmetric bilinear map $\phi_{2}:\g_{0} \times \g_{0} \to
  \g_{1}'$,
\end{itemize}
such that the following equations hold for all $x,y,z \in \g_0$, $u
\in \g_{1}:$
\begin{enumerate}
\item $l_1' (\phi_{2}(x,y)) = \phi_{0}[x,y] - [\phi_{0}(x), \phi_{0}(y)]$,
\item $\phi_{2}(x,l_1(u)) = \phi_{1}[x,u] - [\phi_{0}(x), \phi_{1}(u)]$,
\item $ [\phi_2(x,y), \phi_0(z)]  + \phi_2([x,y],z) + \phi_1(l_3(x,y,z))
= $ $ l_3(\phi_0(x),\phi_0(y), \phi_0(z))  + [\phi_0(x), \phi_2(y,z)] + [\phi_2(x,z), \phi_0(y)] + \phi_2(x, [y,z]) + \phi_2([x,z],y)$.
\end{enumerate}
\label{defi:Lhomo}
\end{defi}
There is a category  {\bf2Term${\bf L}_{\infty}$}  with
$2$-term $L_{\infty}$-algebras as objects and
$L_\infty$-homomorphisms as morphisms.

The Rota-Baxter operators on $L_\infty$-algebras are introduced in \cite{LST}. In the following, we introduction the notion of Rota-Baxter operators on $2$-term $L_\infty$-algebras.
\begin{defn}\label{defi:2RBT}
Let $\huaG=(\g_{1},\g_0,l_1,l_2,l_3)$ be a $2$-term $L_\infty$-algebra. A triple $\frkR=(R_0,R_1,R_2)$, where $R_0:\g_0\longrightarrow \g_0$, $R_1:\g_1\longrightarrow \g_1$ is a chain map, and $R_2:\wedge^2 \g_0\longrightarrow \g_1$ is a linear map, is called a {\bf  Rota-Baxter operator} on $\huaG$ if for all $x,y,x_1,x_2,x_3\in \g_0$ and $u\in \g_1$, the following conditions are satisfied:
\begin{enumerate}
    \item $R_0\big(l_2(R_0 x,y)+l_2(x,R_0 y)\big)-l_2(R_0x,R_0 y)=l_1R_2(x,y)$,
    \label{it:2RBTi}
    \item $R_1\big(l_2(R_1 u,x)+l_2(u, R_0x)\big)-l_2(R_1u,R_0 x)=R_2(l_1(u),x)$,
    \label{it:2RBTii}
    \item $l_2(R_0x_1,R_2(x_2,x_3))+R_2\big(x_3,l_2(R_0x_1,x_2)-l_2(R_0x_2,x_1)\big)
+R_1\big(l_2(R_2(x_2,x_3),x_1)$\\$-l_3(R_0 x_2,R_0x_3,x_1)\big)+c.p.+l_3(R_0x_1,R_0 x_2,R_0 x_3)=0.$
    \label{it:2RBTiii}
  \end{enumerate}
  Moreover, a $2$-term $L_\infty$-algebra $\huaG$ with a triple $\huaR=(R_0,R_1,R_2)$ is
called a {\bf $2$-term Rota-Baxter $L_\infty$-algebra}. We denote a $2$-term Rota-Baxter $L_\infty$-algebra by $(\huaG,\frkR)$. A $2$-term Rota-Baxter $L_\infty$-algebra is called  {\bf strict} if $l_3=0$ and $R_2=0$.
\end{defn}

\begin{defi}\label{defi:RBLhomo}
Let $(\huaG,\frkR)$ and $(\huaG',\frkR')$ be $2$-term Rota-Baxter $L_\infty$-algebras. A {\bf Rota-Baxter $L_\infty$-homomorphism} $\phi=(\phi_0,\phi_1,\phi_2,\phi_3):(\huaG,\frkR)\rar(\huaG',\frkR')$ consists of  a homomorphism $(\phi_0,\phi_1,\phi_2)$ from the $2$-term $L_\infty$-algebra $\huaG$ to the $2$-term $L_\infty$-algebra $\huaG'$ and a linear map $\phi_{3}:\g_{0} \to  \g_{1}'$,
such that, for all $x,y \in \g_0$, $u
\in \g_{1}$, the following equations hold
\begin{eqnarray}
&&l_1' (\phi_{3}(x)) = -R_{0}'(\phi_0(x))+\phi_0(R_0(x));\label{eq:RBLh1}\\
&&\phi_3(l_1(u))=\phi_1(R_1(u))-R_1'(\phi_1(u));\label{eq:RBLh2}\\
&& R_2'(\phi_0(x),\phi_0(y))+R_1'[\phi_3(x),\phi_0(y)]+R_1'[\phi_0(x),\phi_3(y)]\notag\\
&&+R_1'(\phi_2(R_0(x),y))+R_1'(\phi_2(x,R_0(y)))+\phi_3(R_0(x),y)+\phi_3(x,R_0(y))\notag\\
&=&[\phi_3(x),\phi_3(y)]+\phi_2(R_0(x),R_0(y))+\phi_1(R_2(x,y)).\label{eq:RBLh3}
\end{eqnarray}
\end{defi}

It is straightforward to obatain
\begin{prop}
There is a category  {\bf 2TermRB${\bf L_\infty}$}   with
$2$-term Rota-Baxter $L_\infty$-algebras as objects and
Rota-Baxter $L_\infty$-homomorphisms as morphisms.
\end{prop}

\subsection{The equivalence of Rota-Baxter Lie $2$-algebras and $2$-term Rota-Baxter $L_\infty$-algebras }\label{subsubsec:RBLie2TRBL}

The well-known fact between Lie $2$-algebras and $2$-term $L_\infty$-algebras is given as follows.
\begin{theorem}(\cite{BC})\label{thm:equivalent Lie2}
The categories $\LieAAlg$ and $\TTermL$ are equivalent.
\label{them:Lie2TL}
\end{theorem}

In order to prove the main Theorem~\ref{them:RBLie2TRBL} in the part, we recall the constructions of equivalence between $\LieAAlg$ and $\TTermL$ as follows.

The functor from $\LieAAlg$ to $\TTermL$ is denoted as
\begin{equation}
\SLT:\LieAAlg\rar\TTermL.
\label{eq:2LT2}
\end{equation}
Suppose that $L$ is a Lie $2$-algebra. The corresponding $2$-term $L_\infty$-algebra $\SLT(L)=(\g_1,\g_2,l_1,l_2,l_3)$ is given by
\begin{eqnarray*}
&& \g_{0} = L_{0},\quad  \g_{1} = \ker(s) \subseteq L_{1},\\
&&  l_1(u) = t(u)\quad \text{for}\quad u \in \g_1 ,\\
&& l_{2}(x,y) = [x,y]\quad \text{for}\quad x,y \in \g_0 ,\\
&& l_{2}(x,u) = -l_{2}(u,x) = [1_x, u]\quad \text{for}\quad x \in \g_0,u \in \g_1,\\
&& l_2(u,v) = 0\quad \text{for}\quad u,v \in \g_1,\\
&& l_{3}(x,y,z) = \overrightarrow{ {J}_{x,y,z}}\quad \text{for}\quad x,y,z \in \g_0 .
\end{eqnarray*}

Let $L$ and $L'$ be two Lie $2$-algebras. Let $\SLT(L)=(\g_1,\g_2,l_1,l_2,l_3)$ and  $\SLT(L')=(\g_1',\g_2',l_1',l_2',l_3')$ be the corresponding $2$-term $L_\infty$-algebras.  Assume that $F:L \to L'$ is a Lie $2$-algebra homomorphism. The corresponding $L_\infty$-homomorphism $\phi=\SLT(F): \SLT(L) \to \SLT(L')$ is given  by
\begin{eqnarray}
\label{eq:mor1}&& \phi_{0}:\g_0\rar\g_0'\quad \text{by}\quad \phi_{0}(x)= F_{0}(x),\\
\label{eq:mor2}&&  \phi_{1}:\g_1\rar\g_1'\quad \text{by}\quad  \phi_{1}(u)=F_{1}|_{\ker(s)}(u),\\
\label{eq:mor3}&&  \phi_{2}:V_{0} \times V_{0} \to V_{1}'\quad \text{by}\quad  \phi_{2}(x,y)= \vec{F}_{2}(x,y).
\end{eqnarray}

The functor from $\TTermL$ to $\LieAAlg$ is denoted as
\begin{equation}
\TTL:\TTermL\rar\LieAAlg.
\label{eq:2TL2}
\end{equation}
Given a $2$-term $L_\infty$-algebra $\huaG=(\g_1,\g_2,l_1,l_2,l_3)$, we have a Lie $2$-algebra $\TTL(\huaG)=L$, where the object $ L_{0}=\g_0$, the morphism $L_{1}=\g_{0} \oplus \g_{1}$, the source, target, identity-assigning and composite maps are given by
\begin{eqnarray*}
  s(f) &=&  x,  \quad f=(x, \vec{f})\in L_1,\\
  t(f) &=&  x + l_1(\vec{f}),\quad f=(x, \vec{f})\in L_1,\\
  i(y) &=& (y, 0),\quad y\in L_0,\\
  hg   & =& (z, \vec{g} + \vec{h}).
\end{eqnarray*}
Then we see $t(f) - s(f) = l_1(\vec{f})$.
The bracket functor $[\cdot, \cdot]: L \times L \rightarrow L$ is given by
\begin{eqnarray}
  [x,y] &=&  l_2(x,y) , \label{eq:Lb0}\\
{[f,g]} &=& (l_2(x,z), l_2(\vec{f}, z) + l_2(y, \vec{g})) \label{eq:Lb1}\\
        &=& (l_2(x,z), l_2(x,\vec{g}) + l_2(\vec{f},w)) ,\notag
\end{eqnarray}
for arbitrary objects $x,y \in L_0$, and arbitrary morphisms $f:x\rar y,g:z\rar w\in L_1$.
Note that the identity
$$ l_2(\vec{f}, z) + l_2(y, \vec{g}) = l_2(x,\vec{g}) + l_2(\vec{f},z) $$ holds since $l_2(l_1(\vec{f}),\vec{g}) = l_2(\vec{f}, l_1(\vec{g}))$. The Jacobiator for $L$ is given by
$$J_{x,y,z} = ([[x,y],z], l_{3}(x,y,z)).  $$

For each $L_\infty$-homomorphism $\phi:\huaG\rar\huaG'$, we let $\TTL(\huaG)=(\g_1,\g_2,l_1,l_2,l_3)$ and $\TTL(\huaG')=(\g_1',\g_2',l_1',l_2',l_3')$. The Lie $2$-algebra homomorphism $\TTL(\phi)=F:\TTL(\huaG)\rar\TTL(\huaG')$ is defined as follows
\begin{eqnarray*}
  &&F_0 :\g_0\rar \g_0'  \quad\text{by}\quad F_0(x)=\phi_0(x), \\
&&F_1:\g_1\rar \g_1' \quad\text{by}\quad F_1(f)=F_1(x,\vec{f})=(\phi_0(x),\phi_1(\vec{f})),\\
&&F_{2}:\g_0\times\g_0\rar \g_1' \quad\text{by}\quad F_{2}(x,y) = ([\phi_{0}(x), \phi_{0}(y)], \phi_{2}(x,y)) .
\end{eqnarray*}

Finally, the natural isomophisms $\alpha:\TTL\SLT \Longrightarrow 1_{\LieAAlg}$ and
$\beta:\SLT\TTL \Longrightarrow 1_{\TTermL}$ imply the equivalence between $\LieAAlg$ and $\TTermL$.

As a generalization of Theorem~\ref{thm:equivalent Lie2}, we have
\begin{theorem}
The categories $\RBLieAAlg$ an $\TTermRBL$ are equivalent.
\label{them:RBLie2TRBL}
\end{theorem}
\begin{proof}
First we construct a functor $\SLT^{\rm RB}:\RBLieAAlg\rar\TTermRBL$ which `lifts' the functor $\SLT$ in~$($\ref{eq:2LT2}$)$ in form of the following commutative diagram
\[\xymatrix{
\RBLieAAlg\ar[r]^{\SLT^{\rm RB}}\ar[d]_{U_\RBLieAAlg}& \TTermRBL\ar[d]^{U_\TTermRBL}\\
\LieAAlg\ar[r]_{\SLT}& \TTermL
}\]
where $U_\RBLieAAlg$ and $U_\TTermRBL$ are forgetful functors.

Given a Rota-Baxter Lie $2$-algebra $(L,P)$, we obtain a $2$-term Rota-Baxter $L_\infty$ algebra $\SLT^{\rm RB}(L,P)=(\huaG, \frkR)$. Here $\huaG=(\g_1,\g_0,l_1,l_2,l_3)$ is the $2$-term $L_\infty$ algebra $\SLT(P)$, and $\huaR=(R_0, R_1, R_2)$ on $\huaG$ is given by %defined as follows. For all $x, y\in\g_0, h\in\g_1,$
\begin{eqnarray*}
&&R_0 :\g_0\rar \g_0  \quad\text{by}\quad R_0(x)=P_0(x), \\
&&R_1:\g_1\rar \g_1 \quad\text{by}\quad R_1(u)=P_1(u),\\
&&R_{2}:\g_0\times\g_0\rar \g_1' \quad\text{by}\quad R_{2}(x,y) =\overrightarrow{\huaR}_{x,y}.
\end{eqnarray*}%Item~(\mref{it:rb'2})
In the following, we first show that the conditions~(\ref{it:2RBTi}), (\ref{it:2RBTii}) and (\ref{it:2RBTiii}) in Definition~\ref{defi:2RBT} hold.

The condition ~$($\ref{it:2RBTi}$)$ holds since
$$R_0\left(l_2(R_0 x,y)+l_2(x,R_0 y)\right)-l_2(R_0x,R_0 y)=(t-s)\huaR_{x,y}=t\overrightarrow{\huaR_{x,y}}=l_1R_2(x,y).$$

The naturality of $\huaR_{x,y}$ implies that for any $f:x\rar z$, we have the identity
\begin{equation}
\huaR_{x,y}\left(P_1[P_1(f),1_y]+P_1[f,1_{P_0(y)}]\right)=[P_1(f),1_{P_0(y)}] \huaR_{z,y},
\label{eq:ntRB}
\end{equation}
Taking the arrow parts of both sides of the above Eq.~$($\ref{eq:ntRB}$)$, we have
\begin{equation*}
\overrightarrow{\huaR_{x,y}}+\left(\overrightarrow{P_1[P_1(f),1_y]}
+\overrightarrow{P_1[f,1_{P_0(y)}]}\right)
=\overrightarrow{[P_1(f),1_{P_0(y)}]}+\overrightarrow{\huaR_{z,y}},
%\mlabel{eq:ntRB1}
\end{equation*}
which implies that
\begin{equation}
P_1\left([P_1(\vec{f}),1_y]+ [\vec{f},1_{P_0(y)}]\right)-[P_1(\vec{f}),1_{P_0(y)}]
=\overrightarrow{\huaR_{z-x,y}}.
\label{eq:ntRB1}
\end{equation}
Thus we have
\begin{equation}
R_1\left(l_2(R_1 (\vec{f}),y)+l_2(\vec{f}, R_0(y))\right)-l_2(R_1(\vec{f}),R_0 (y))
=R_2(l_1(\vec{f}),y).
\end{equation}
This implies that the conditon~$($\ref{it:2RBTii}$)$ holds.

It is straightforward to check that Eq.$($\ref{eq:RBcoh}$)$ is equivalent to
\begin{eqnarray*}
  [1_{P_{0}(x)},\overrightarrow{\huaR_{y,z}}]+ \overrightarrow{\huaR_{z,[P_{0}(x),y]+[x,P_{0}(y)]}}   +P_{1}\left([\overrightarrow{\huaR_{y,z}},1_{x}]-\overrightarrow{J_{P_0(y),P_0(z),x}}\right)
+c.p.+\overrightarrow{J_{P_{0}(x),P_{0}(y),P_{0}(z)}}=0,
\end{eqnarray*}
which implies that
\begin{eqnarray*}
    &&l_2(R_0x,R_2(y,z))+R_2\left(g(z,l_2(R_0x,y)-l_2(R_0y,x )\right)\\
    &&+R_1\left(l_2(R_2(y,z),x)-l_3(R_0 y,R_0z,x )\right)+c.p.+l_3(R_0x,R_0 y,R_0 z)=0.
\end{eqnarray*}
This implies that the conditon~$($\ref{it:2RBTiii}$)$ holds.

Next we construct a Rota-Baxter $L_\infty$-homomorphism from a Rota-Baxter Lie $2$-algebra homomorphism. Let $F:(L,P,\huaR)\rar (L',P',\huaR')$ be a Rota-Baxter Lie $2$-algebra homomorphism. Let $(\huaG,\frkR)=\SLT^{\rm RB}(L,P,\huaR)$ and $(\huaG,\frkR)=\SLT^{\rm RB}(L',P',\huaR')$. Then we obtain an $L_\infty$-homomorphism $\phi=\SLT(F):\huaG\rar\huaG'$ of $\SLT^{\rm RB}(F)$ as in ~(\ref{eq:mor1})-\eqref{eq:mor3}. Define a map $\phi_3:V_0\rar V_1'$ by $$\phi_3(x)=\overrightarrow{F_3(x)}:0\rar -P_{0}'(\phi_0(x))+\phi_0(P_0(x)).$$

 Eq.~$($\ref{eq:RBLh1}$)$ holds in Definition~\ref{defi:RBLhomo} since
$$l_1' (\phi_{3}(x)) =t(\overrightarrow{F_3(x)})= -P_{0}'(\phi_0(x))+\phi_0(P_0(x))= -R_{0}'(\phi_0(x))+\phi_0(R_0(x)).$$

By the naturality of $F_3$, for every morphism $f:x\rar y $, we have a $$\phi_1(P_1(f))F_3(x)=F_3(y)P_1'(\phi_1(f)).$$
Furthermore, we have
\begin{eqnarray*}
\overrightarrow{F_3(x)}+\phi_1(P_1(\overrightarrow{f}))=\overrightarrow{\phi_1(P_1(f))F_3(x)}
=\overrightarrow{F_3(y)P_1'(\phi_1(f))}
=P_1'(\phi_1(\overrightarrow{f}))+\overrightarrow{F_3(y)},
\end{eqnarray*}
which implies that $$\phi_3(l_1(f))=\overrightarrow{F_3(y-x)}=\overrightarrow{F_3(y)}-\overrightarrow{F_3(x)}
=\phi_1(P_1(\overrightarrow{f}))-P_1'(\phi_1(\overrightarrow{f})).$$
Thus for  any $u\in \g_1$, we have
$$\phi_3(l_1(u))= \phi_1(P_1(u))-P_1'(\phi_1(u))= \phi_1(R_1(u))-R_1'(\phi_1(u)).$$
This implies that Eq.~$($\ref{eq:RBLh2}$)$ in Definition~\ref{defi:RBLhomo} holds.

It is straightforward to check that Eq.~$($\ref{eq:RBcohm}$)$ is equivalent to the following equation
\begin{eqnarray*}
&&(\overrightarrow{F_3[P_0(x),y]}+\overrightarrow{F_3[x,P_0(y)]})
(P_1'\overrightarrow{F_2(P_0(x),y)} +P_1'\overrightarrow{F_2(x,P_0(y))})\notag\\
&&(P_1'[\overrightarrow{F_3(x)},1_{F_0(y)}]  +P_1'[1_{F_0(x)},\overrightarrow{F_3(y)}])\overrightarrow{RB_{F_0(x),F_0(y)}}\notag\\
&=&F_1(\overrightarrow{RB_{x,y}})\overrightarrow{F_2(P_0(x),P_0(y))}
[\overrightarrow{F_3(x)},\overrightarrow{F_3(y)}],\label{eq:RBcohmarr}
\end{eqnarray*}
which implies that Eq.~$($\ref{eq:RBLh3}$)$ in Definition~\ref{defi:RBLhomo} holds.

One can also deduce that $\SLT^{\rm RB}$ preserves the identity homomorphisms and the composition of homomorphisms. Thus $\SLT^{\rm RB}$ is a functor from $\RBLieAAlg$ to $\TTermRBL$.

Conversely, we construct a functor $\TTL^{\rm RB}:\TTermRBL\rar\RBLieAAlg$ as a `lifting' of the functor $\TTL$ in~$($\ref{eq:2TL2}$)$ in form of the following commutative diagram
\[\xymatrix{
\TTermRBL\ar[r]^{\TTL^{\rm RB}}\ar[d]_{U_\TTermRBL}&\RBLieAAlg \ar[d]^{U_\RBLieAAlg}\\
\TTermL\ar[r]_{\TTL}& \LieAAlg
}\]
where $U_\TTermRBL$ and $U_\RBLieAAlg$ are corresponding forgetful functors.

Let $(\huaG,\frkR)$ be a $2$-term Rota-Baxter $L_\infty$ algebra, where $\huaG=(\g_0, \g_1, l_1, l_2, l_3)$ is a $2$-term $L_\infty$ algebra and $\frkR=(R_0, R_1, R_2)$ is a Rota-Baxter operator on $\huaG$. Then we have a  Lie $2$-algebra $\TTL(\huaG)=L$ with $L_0=\g_0$ and $L_1=\g_0\oplus\g_1$. Define a linear functor $P:L\rar L$ by
\begin{eqnarray*}
&&P_0:L_0\rar L_0,\quad P_0(x)=R_0(x)\, \quad \forall~x\in L_0,\\
&&P_1:L_1\rar L_1,\quad P_1(y,u)=(R_0(y), R_1(u)),\quad \forall~y\in \g_0, u\in \g_1.
\end{eqnarray*}
The natural isomorphism $\huaR_{x, y}:[P_{0}(x),P_{0}(y)]\rar P_{0}([P_{0}(x),y]+[x,P_{0}(y)])$ is defined by
 $$\huaR_{x,y}=([P_0(x),P_0(y)], R_2(x,y)).$$
By the various conditions of $2$-term Rota-Baxter $L_\infty$ algebra, we can deduce that $(L,P,\huaR)$ is a Rota-Baxter Lie $2$-algebra. Thus, we have constructed a Rota-Baxter Lie $2$-algebra $(L,P,\huaR)=\TTL^{\rm RB}(\huaG,\frkR)$ from a $2$-term Rota-Baxter $L_\infty$ algebra $(\huaG,\frkR)$.

For any Rota-Baxter $L_\infty$-homomorphism $\phi:(\huaG,\huaR)\rar (\huaG',\huaR')$, next we construct a Rota-Baxter Lie $2$-algebra homomorphism $F=\TTL^{\rm RB}(\phi)$  from $\TTL^{\rm RB}(\huaG,\frkR)$ to $\TTL^{\rm RB}(\huaG',\frkR')$.

The underlying Lie $2$-algebra homomorphism is given by
\begin{eqnarray*}
&&F_0=\phi_0:L_0\rar L_0',\\
&&F_1=\phi_0\oplus\phi_1:L_1\rar L_1',\\
&&F_2:L_0\times L_0\rar L_1',\quad F_2(x,y)=([\phi_0(x),\phi_0(y)],\phi_2(x,y)).
\end{eqnarray*}
The natural transformation $F_3(x):P_0'(F_0(x))\rar F_0(P_0(x))$ is defined by
$$F_3(x)=(P_0'(F_0(x)), \phi_3(x)).$$
Applying the correspondence between the composition of morphisms and the addition of their arrow parts, Eq.~$($\ref{eq:RBLh2}$)$ implies the naturality of $F_3$. Eq.~$($\ref{eq:RBcohm}$)$ also holds by Eq.~$($\ref{eq:RBLh3}$)$.
 Thus $F$ is a Rota-Baxter Lie $2$-algebra homomorphism. Furthermore, $\TTL^{\rm RB}$ preserves the identity homomorphisms and the composition of homomorphisms. Therefore, $\TTL^{\rm RB}$ is a functor from $\TTermRBL$ to $\RBLieAAlg$.

We are left to show that there are natural isomorphisms $$\alpha^{\rm RB}:\TTL^{\rm RB}\SLT^{\rm RB}\Rightarrow 1_{\RBLieAAlg}\quad\text{and}\quad \beta^{\rm RB}:\SLT^{\rm RB}\TTL^{\rm RB}\Rightarrow 1_{\TTermRBL}.$$

For any Rota-Baxter Lie $2$-algebra $(L,P,\huaR)$, we obtain a $2$-term Rota-Baxter $L_\infty$ algebra  $$\SLT^{\rm RB}(L,P,\huaR)=(\huaG,\frkR)=((\g_0,\g_1,l_1,l_2,l_3),(R_0,R_1,R_2)),$$
 where $\SLT(L)=\huaG$, and $R_0=P_0$, $R_1=P_1|_{\g_1}$, $R_2(x,y)=\overrightarrow{\huaR_{x,y}}$. Applying the functor $\TTL^{\rm RB}$ to $(\huaG,\frkR)$, we obtain a  Rota-Baxter Lie $2$-algebra, denoted by $(L',P',\huaR')$. Here $L'= \TTL(\SLT(L))$, and for all $x\in L_0'$ and $(y, u)\in L_1'=\g_0\oplus\g_1$, one has
 $$P_0'(x)=R_0(x), P_1'(y,u)=(R_0(y),R_1(u)), \huaR_{x,y}'=([P_0(x),P_0(y)],R_2(x,y))=\huaR_{x,y}.$$ By the isomorphism $\alpha_L:L'\rar L$ of semistrict Lie $2$-algebras: $(\alpha_L)_0(x)=x$ and $(\alpha_L)_1(y,u)=i(y)+u$, we have $P_0'(x)=P_0(x)$, and $$P_1((\alpha_L)_1(y,u))=P_1(i(y)+u)=i(P_0(y))+P_1(u)=(\alpha_L)_1(P_0(y),P_1(u))
=(\alpha_L)_1(P_1'(y,u)).$$
Thus $\alpha^{\rm RB}:(L',P',\huaR')\rar (L,P,\huaR)$ is an isomorphism of Rota-Baxter Lie $2$-algebras. Also by the naturality of $\alpha$, we see that $\alpha^{\rm RB}$ is a natural isomorphism.

For a $2$-term Rota-Baxter $L_\infty$ algebra $(\huaG,\huaR)=((\g_0,\g_1,l_1,l_2,l_3),(R_0,R_1,R_2))$, applying the functor $\TTL^{\rm RB}$ to $(\huaG,\huaR)$, we obatin a Rota-Baxter Lie $2$-algebra $(L,P,\huaR)$, where $L_0=\g_0$, $L_1=\g_0\oplus\g_1$, $P_0(x)=R_0(x)$, $P_1(y,u)=R_0(y)+R_1(u)$ and  $\huaR_{x,y}=([P_0(x),P_0(y)],R_2(x,y))$ for all $x\in L_0$ and $(y, u)\in \g_0\oplus\g_1$. Applying $\SLT^{\rm RB}$ to $(L,P)$, we have a $2$-term Rota-Baxter $L_\infty$ algebra $(\huaG',\frkR')$, where $\huaG'=\SLT(\TTL(\huaG))$, $R_0'(x)=P_0(x)$, $R_1'(u)=P_1(u)=R_1(h)$ and $R'_2(x,y)=\overrightarrow{\huaR_{x,y}}=R_2(x,y)$ for any $x,y\in\g_0$ and $u\in\g_1$.  Thus we obtain the $2$-term Rota-Baxter $L_\infty$-algebra isomorphism $\beta^{\rm RB}:(\huaG',\huaR')\rar(\huaG,\huaR)$. The naturality of $\beta^{\rm RB}$ follows that of $\beta$. Then we obtain a natural isomorphism $\beta^{\rm RB}$.
\end{proof}

\begin{remark}
We can further obtain $2$-categories $\RBLieAAlg$ and $\TTermRBL$ by introducing $2$-morphisms and strengthen Theorem~\ref{them:RBLie2TRBL} to the $2$-equivalence of $2$-categories.
\end{remark}

For strict Rota-Baxter Lie $2$-algebras, there is a category $\BSRBLieAAlg$ with strict Rota-Baxter Lie $2$-algebra as objects and Rota-Baxter Lie $2$-algebra homomorphisms as morphisms, which is a subcategory of $\RBLieAAlg$.

 For strict $2$-term $L_{\infty}$-algebra,  there is a category $\SRBTTermL$  with
strict $2$-term $L_{\infty}$-algebras as objects and
$L_\infty$-homomorphisms as morphisms, which is a subcategory of $\TTermRBL$.

It is straightforward to check that
\begin{cor}
The categories SRBLie2Alg and SRB2Term$L_\infty$ are equivalent.
\label{them:Lie2TL}
\end{cor}

\section{Strict $2$-term Rota-Baxter $L_{\infty}$-algebras and crossed modules of Rota-Baxter Lie algebras}\label{sec:sRBLie}
 In this section, we study the relations between strict $2$-term Rota-Baxter $L_{\infty}$-algebras and  crossed modules of Rota-Baxter Lie algebras.

  First, we recall the definition of crossed modules of Lie algebras.
\begin{defi}
  A {\bf crossed module of Lie algebras} is a quadruple $((\g_1,[\cdot,\cdot]_{\g_1}),(\g_0,[\cdot,\cdot]_{\g_0}),\dM,\rho)$, where $(\g_1,[\cdot,\cdot]_{\g_1})$ and $(\g_0,[\cdot,\cdot]_{\g_0})$ are Lie algebras, $\dM:\g_1\rightarrow \g_0$ is a Lie algebra homomorphism and $\rho:\g_0\rightarrow\Der(\g_1)$ is an action of Lie algebra $\g_0$ on Lie algebra $\g_1$ as a derivation, such that
  \begin{eqnarray}
    \dM(\rho(x)(u))=[x,\dM u]_{\g_0},\quad \rho(\dM u)(v)=[u,v]_{\g_1},\quad \forall~x\in\g_0,u,v\in\g_1.
  \end{eqnarray}
\end{defi}

\emptycomment{\begin{ex}
  Let  $(\g,[\cdot,\cdot]_\g)$ be a Lie algebra and $\h$ an ideal of $\g$. Then $(\g,\h,\dM=\imath,\rho=\ad)$ is a crossed module for Lie algebras, where $\imath:\h\longrightarrow \g$ is the inclusion.
\end{ex}

 \begin{ex}
  For any Lie algebra homomorphism $f:\g\longrightarrow \h$, $\Ker~ f$ is an ideal of $\g$. Then $(\g,\Ker ~f,\imath,\ad)$ is a crossed module of Lie algebras.
\end{ex}

\begin{pro}
  Let  $((\g_0,[\cdot,\cdot]_{\g_0}),(\g_1,[\cdot,\cdot]_{\g_1}),\dM,\rho)$ be a crossed module of Lie algebras. Then there is a Lie algebra structure on $\g_0\oplus\g_1$ given by
  \begin{equation}
    [x+u,y+v]=[x,y]_{\g_0}+\rho(x)v-\rho(y)u+[u,v]_{\g_1},\quad\forall~x,y\in\g_0,u,v\in\g_1.
  \end{equation}
\end{pro}}

Relations between strict Lie $2$-algebras and crossed modules of Lie algebras are described in the following theorem.
\begin{thm}{\rm (\cite{BC})}
  There is a one-to-one corresponding between strict Lie $2$-algebras and crossed modules of Lie algebras.
\end{thm}

\begin{defi}
  A {\bf crossed module of Rota-Baxter Lie algebras} is a quadruple \\ $((\g_1,[\cdot,\cdot]_{\g_1},T_1),(\g_0,[\cdot,\cdot]_{\g_0},T_0),\dM,\rho)$, where $(\g_1,[\cdot,\cdot]_{\g_1},T_1)$ and $(\g_0,[\cdot,\cdot]_{\g_0},T_0)$ are Rota-Baxter Lie algebras, $\dM:\g_1\rightarrow \g_0$ is a Rota-Baxter Lie algebra homomorphism and $(\rho,T_1):\g_0\rightarrow\Der(\g_1)$ is an action of Rota-Baxter Lie algebra $(\g_0,T_0)$ on Lie algebra $\g_1$ as a derivation of the Lie algebra, such that
  \begin{eqnarray}
    \dM(\rho(x)(u))=[x,\dM u]_{\g_0},\quad \rho(\dM u)(v)=[u,v]_{\g_1},\quad \forall~x\in\g_0,u,v\in\g_1.
  \end{eqnarray}
\end{defi}
It is obvious that $((\g_0,[\cdot,\cdot]_{\g_0}),(\g_1,[\cdot,\cdot]_{\g_1}),\dM,\rho)$ is a crossed module of Lie algebras.

\begin{ex}
  Let  $(\g,[\cdot,\cdot]_\g,R)$ be a Rota-Baxter Lie algebra and $\h$ a Rota-Baxter Lie ideal of $(\g,R)$. Then $(\g,\h,\dM=\imath,\rho=\ad)$ is a crossed module for Rota-Baxter Lie algebras, where $\imath:\h\longrightarrow \g$ is the inclusion.
\end{ex}

\begin{ex}
  For any Rota-Baxter Lie algebra homomorphism $f:\g\longrightarrow \h$,  $(\g,\Ker ~f,\imath,\ad)$ is a crossed module of Rota-Baxter Lie algebras.
\end{ex}

\begin{pro}
  Let  $((\g_0,[\cdot,\cdot]_{\g_0},T_0),(\g_1,[\cdot,\cdot]_{\g_1},T_1),\dM,\rho)$ be a crossed module of Rota-Baxter Lie algebras. Then there is a Rota-Baxter Lie algebra structure on $\g_0\oplus\g_1$ given by
  \begin{eqnarray}
    {[x+u,y+v]}&=&[x,y]_{\g_0}+\rho(x)v-\rho(y)u+[u,v]_{\g_1},\\
    T(x+u)&=&T_0(x)+T_1(u),\quad\forall~x,y\in\g_0,u,v\in\g_1.
  \end{eqnarray}
\end{pro}
\begin{proof}
  By the fact that $((\g_0,[\cdot,\cdot]_{\g_0}),(\g_1,[\cdot,\cdot]_{\g_1}),\dM,\rho)$ is a crossed module of Lie algebras, $(\g_0\oplus\g_1,[\cdot,\cdot])$ is a Lie algebra.

  Furthermore, it is straightforward to check that $T$ is a Rota-Baxter operator on the Lie algebra $\g\oplus\h$ if and only if $T_0$ is a Rota-Baxter operator on the Lie algebra $\g_0$, $T_1$ is a Rota-Baxter operator on the Lie algebra $\g_1$ and the following equation holds:
  $$T_1(\rho(T_0 x)u+\rho(x)T_1u)=\rho(T_0 x)T_1 u,\quad \forall~x\in\g_0,u\in \g_1,$$
  which follows from that $(\rho,T_1)$ is a representation of the Rota-Baxter Lie algebra $(\g_0,T_0)$ on $\g_1$.
\end{proof}

\begin{thm}
  There is a one-to-one corresponding between strict $2$-term Rota-Baxter $L_\infty$-algebras and crossed modules of Rota-Baxter Lie algebras.
\end{thm}
\begin{proof}
  Let $(\g_0,\g_1,l_1,l_2,l_3=0;R_0,R_1,R_2=0)$ be a strict $2$-term Rota-Baxter $L_\infty$-algebra. Define the brackets $[\cdot,\cdot]_{\g_0}$ and  $[\cdot,\cdot]_{\g_1}$ by
  \begin{equation*}
    [x,y]_{\g_0}=l_2(x,y),\quad [u,v]_{\g_1}=l_2(l_1(u),v),\quad\forall~x,y\in\g_0,u,v\in \g_1.
  \end{equation*}
Define $\rho:\g_0\rightarrow \gl(\g_1)$ by
\begin{equation*}
  \rho(x)u=l_2(x,u),\quad\forall~x\in\g_0,u\in\g_1.
\end{equation*}
Then $((\g_0,[\cdot,\cdot]_{\g_0}),(\g_1,[\cdot,\cdot]_{\g_1}),\dM=l_1,\rho)$ is a a crossed module of Lie algebras.

Set $T_0=R_0$ and $T_1=R_1$. By condition (a) in Definition \ref{defi:2RBT}, $T_0$ is a Rota-Baxter operator on the Lie algebra $\g_0$. By condition (b) in Definition \ref{defi:2RBT} and condition (a) in Definition \ref{defi:2-term}, for $u,v\in\g_1$, we have
\begin{eqnarray*}
 && T_1([T_1(u),v]_{\g_1}+[u,T_1 (v)]_{\g_0})-[T_1(u),T_1(v)]_{\g_1}\\
 &=&R_1\big(l_2(l_1 R_1(u),v)+l_2(l_1(u),R_1(v)\big)-l_2(l_1 R_1(u),R_1(v)))\\
 &=&R_1\big(l_2(R_1(u),l_1 (v))+l_2(u,R_0l_1(v)\big)-l_2( R_1(u),R_0l_1(v)))\\
 &=&0,
\end{eqnarray*}
which implies that $T_1$ is a Rota-Baxter operator on the Lie algebra $\g_1$. By the fact that $\dM$ is a Lie algebra homomorphism from $\g_1$ to $\g_0$ and $l_1\circ R_1=R_0\circ l_1$, $\dM$ is a Rota-Baxter Lie algebra homomorphism from $(\g_1,T_1)$ to $(\g_0,T_0)$. By condition (b) in Definition \ref{defi:2RBT}, $(\rho,T_1):\g_0\rightarrow\Der(\g_1)$ is an action of Rota-Baxter Lie algebra $(\g_0,T_0)$ on Lie algebra $\g_1$. Therefore,
$((\g_0,[\cdot,\cdot]_{\g_0},T_0),(\g_1,[\cdot,\cdot]_{\g_1},T_1),\dM=l_1,\rho)$ is a crossed module of Rota-Baxter Lie algebras.

Conversely, a crossed module of Rota-Baxter Lie algebras $((\g_0,[\cdot,\cdot]_{\g_0},T_0),(\g_1,[\cdot,\cdot]_{\g_1},T_1),\dM,\rho)$ gives rise to a strict $2$-term Rota-Baxter $L_\infty$-algebra $(\g_0,\g_1,l_1=\dM,l_2,l_3=0;R_0=T_0,R_1=T_1,R_2=0)$, where $l_2:\g_i\wedge \g_j\rightarrow \g_{i+j},~0\leq i+j\leq1$ is given by
$$l_2(x,y)=[x,y]_{\g_0},\quad l_2(x,u)=\rho(x)u,\quad \forall~x,y\in\g_0,u,v\in\g_1$$
The conditions in crossed module of Rota-Baxter Lie algebras  give various conditions for a strict $2$-term Rota-Baxter $L_\infty$-algebra. We omit the details.
\end{proof}

Let $(\g,\ast)$ be a pre-Lie algebra and $V$  a vector
space. A {\bf representation} of $\g$ on $V$ consists of a pair
$(l,r)$, where $l:\g\longrightarrow \gl(V)$ is a representation
of the Lie algebra $\g^c$ on $V $ and $r:\g\longrightarrow \gl(V)$ is a linear
map satisfying \begin{eqnarray}\label{representation condition 2}
 r_xl_y-r_yl_x=r_{x\ast y}-r_yr_x, \quad \forall~x,y\in \g.
\end{eqnarray}

Recall that a {\bf crossed module of pre-Lie algebras} is a quadruple $((\g_0,\ast_0),(\g_1,\ast_1),\delta,(l,r))$, where $(\g_0,\ast_0)$ and $(\g_1,\ast_1)$ are pre-Lie algebras, $\delta:\g_1\rightarrow \g_0$ is a homomorphism of pre-Lie algebras, and $(l,r)$ is a representation of the pre-Lie algebra $(\g_0,\ast_0)$ on $\g_1$, such that for $x\in\g_0$ and $u,v\in\g_1$ ,the following equalities are satisfied:
\begin{eqnarray}
  \delta(l_x u)&=&x\ast_0 \delta u,\quad \delta(r_x u)=(\delta u)\ast_0 x,\\
  l_{\delta u}v&=&r_{\delta v}u=u\ast_1 v.
\end{eqnarray}

\begin{pro}\label{pro:relation crossed module}
Let  $((\g_0,[\cdot,\cdot]_{\g_0},T_0),(\g_1,[\cdot,\cdot]_{\g_1},T_1),\dM,\rho)$ be a crossed module of Rota-Baxter Lie algebras. Define $\ast_0:\g_0\otimes\g_0\rightarrow \g_0$, $\ast_1:\g_1\otimes\g_1\rightarrow \g_1$ and $l,r:\g_0\rightarrow \gl(\g_1)$ by
\begin{eqnarray*}
  x\ast_0 y&=&[T_0x,y]_{\g_0},\quad u\ast_1v=[T_1 u,v]_{\g_1},\\
  l_x u&=&\rho(T_0 x)u,\quad r_x u=-\rho(x)T_1 (u),\quad \forall~x,y\in\g_0,u,v\in \g_1.
\end{eqnarray*}
Then $((\g_0,\ast_0),(\g_1,\ast_1),\dM,(l,r))$ is a crossed module of pre-Lie algebras.
\end{pro}
\begin{proof}
  Since $T_0$ is a Rota-Baxter operator on the Lie algebra $(\g_0,[\cdot,\cdot]_{\g_0})$, $(\g_0,\ast_0)$ is a pre-Lie algebra. Similarly, $(\g_1,\ast_1)$ is also a pre-Lie algebra. By the fact that $\dM$ is a Rota-Baxter Lie algebra homomorphism, we have
  \begin{eqnarray*}
    \dM(u\ast_1 v)&=&\dM [T_1 u,v]_{\g_1}=[\dM (T_1 u),\dM v]_{\g_1}\\
    &=&[T_0(\dM u),\dM ]_{\g_1}=(\dM u)\ast_0 (\dM v),
  \end{eqnarray*}
  which implies that $\dM$ is a pre-Lie algebra homomorphism from $\g_1$ to $\g_0$.

  By the fact that $\rho$ is a representation of the Lie algebra $\g_0$ on $\g_1$ and $T_0$ is a Rota-Baxter operator on $\g_0$, we have
  \begin{eqnarray*}
    l_{[x,y]_{T_0}}&=&\rho(T_0([x,y]_{T_0}))=\rho([T_0x,T_0y]_{\g_0})\\
    &=&[\rho(T_0 x),\rho(T_0 y)]=[l_x,l_y],
  \end{eqnarray*}
  which implies that $l$ is a representation of the sub-adjacent Lie algebra $\g_0^c$ on $\g_1$. Furthermore, by Eq. \eqref{eq:rep-RB} in the representation of the Rota-Baxter Lie algebra, we have
  \begin{eqnarray*}
  &&  l_x(r_y u)-r_y(l_xu)-r_{x\ast_0 y}u+r_y(r_x u)\\
  &=&-\rho(T_0 x)\rho(y)(T_1 u)+\rho(y)T_1(\rho(T_0 x)u)+\rho([T_0 x,y]_{\g_0})u+\rho(y)T_1\rho(x)(T_1 u)\\
  &=&\rho(y)T_1(\rho(T_0 x)u)-\rho(y)\rho(T_0 x)u+\rho(y)T_1\rho(x)(T_1 u)=0.
  \end{eqnarray*}
  Thus $(l,r)$ is a representation of the pre-Lie algebra $(\g_0,\ast_0)$ on $\g_1$.

  Furthermore, the condition $\dM(\rho(x)(u))=[x,\dM u]_{\g_0}$ implies that
  $$\dM(l_x u)=x\ast_0 \dM u,\quad \dM(r_x u)=(\dM u)\ast_0 x$$
  hold and the condition $\rho(\dM u)(v)=[u,v]_{\g_1}$ implies that
  $$l_{\dM u}v=r_{\dM v}u=u\ast_1 v$$
  hold.
Therefore we obtain a crossed module of pre-Lie algebras $((\g_0,\ast_0),(\g_1,\ast_1),\dM,(l,r))$.
\end{proof}

\begin{pro}\label{pro:pre-Lie crossed module}{\rm(\cite{Sh})}
Let $((\g_0,\ast_0),(\g_1,\ast_1),\dM,(l,r))$ be a crossed module of pre-Lie algebras. Then $((\g_0,[\cdot,\cdot]_{\g_0}),(\g_1,[\cdot,\cdot]_{\g_1}),\dM,\rho=l-r)$ is a crossed module of Lie algebras, where the brackets $[\cdot,\cdot]_{\g_0}$ and $[\cdot,\cdot]_{\g_1}$ are given by
\begin{equation}
  [x,y]_{\g_0}=x\ast_{0} y-y\ast_{0} x,\quad [u,v]_{\g_1}=u\ast_{1} v-v\ast_{1} u
\end{equation}
for $x,y\in\g_0,u,v\in\g_1$.
\end{pro}

Let $((\g_1,[\cdot,\cdot]_{\g_1}),(\g_0,[\cdot,\cdot]_{\g_0}),\dM_\g,\rho_\g)$ and $((\h_1,[\cdot,\cdot]_{\h_1}),(\h_0,[\cdot,\cdot]_{\h_0}),\dM_\h,\rho_\h)$ be two crossed modules of Lie algebras.   Recall that a {\bf homomorphism} from  $(\g_0,\g_1,\dM_\g,\rho_\g)$ to $(\h_0,\h_1,\dM_\h,\rho_\h)$ is a pair $(\psi_0,\psi_1)$, such that $\psi_0:\g_0\rightarrow \h_0$ is a Lie algebra homomorphism and $\psi_1:\g_1\rightarrow \h_1$ is a Lie algebra homomorphism satisfying
\begin{eqnarray}
  \dM_\h\circ \psi_1=\psi_0\circ \dM_\g,\quad \psi_1(\rho_\g(x)v)=\rho_\h(\psi_0(x))\psi_1(v),\quad \forall~x\in\g_0,v\in \g_1.
\end{eqnarray}

\begin{pro}
  Let  $((\g_0,[\cdot,\cdot]_{\g_0},T_0),(\g_1,[\cdot,\cdot]_{\g_1},T_1),\dM,\rho)$ be a crossed module of Rota-Baxter Lie algebras. Then $((\g_0,[\cdot,\cdot]_{T_0}),(\g_1,[\cdot,\cdot]_{T_1}),\dM,\rho_T)$ is a crossed module of Lie algebras, where $[\cdot,\cdot]_{T_0}$, $[\cdot,\cdot]_{T_1}$ and $\rho_T$ are given by
\begin{eqnarray*}
  [x,y]_{T_0}&=&[T_0x,y]_{\g_0}-[T_0y,x]_{\g_0},\\
{[u,v]_{T_1}}&=&[T_1 u,v]_{\g_1}-[T_1 v,u]_{\g_1},\\
  \rho_T(x)u&=&\rho(T_0 x)u+\rho(x)T_1 (u), \quad \forall~x,y\in\g_0,u,v\in\g_1.
\end{eqnarray*}
Furthermore, the pair $(T_0,T_1)$ is a  homomorphism from crossed module of Lie algebras $((\g_0,[\cdot,\cdot]_{T_0}),(\g_1,[\cdot,\cdot]_{T_1}),\dM,\rho_T)$ to  $((\g_0,[\cdot,\cdot]_{\g_0},T_0),(\g_1,[\cdot,\cdot]_{\g_1},T_1),\dM,\rho)$.
\end{pro}
\begin{proof}
The first conclusion follows from Proposition \ref{pro:relation crossed module} and \ref{pro:pre-Lie crossed module}.

Since $T_0$ is a Rota-Baxter operator on $\g_0$, $T_0$ is a Lie algebra homomorphism from $(\g_0,[\cdot,\cdot]_{T_0})$ to $(\g_0,[\cdot,\cdot]_{\g_0})$. Similarly, $T_1$ is a Lie algebra homomorphism from $(\g_1,[\cdot,\cdot]_{T_1})$ to $(\g_1,[\cdot,\cdot]_{\g_1})$. By the fact that $\dM$ is a Rota-Baxter Lie algebra homomorphism, we have $\dM\circ T_1=T_0\circ \dM$ and furthermore, by the fact that $\rho$ is a representation of the Rota-Baxter Lie algebra $\g_0$ on $\g_1$, we have
$$T_1(\rho_T(x)u)-\rho(T_0(x))T_1(u)=T_1(\rho(T_0 x)u+\rho(x)T_1 (u))-\rho(T_0(x))T_1(u)=0.$$
Thus the second conclusion follows.
\end{proof}

{\bf Acknowledgements}: This research was supported by the National Natural Science Foundation of China (11901501), the China Postdoctoral Science Foundation (2021M700750) and the National Key Research and Development Program of China (2021YFA1002000).


\begin{thebibliography}{abcdsfgh}
\bibitem{Ba} G. Baxter, An analytic problem whose solution follows from a simple algebraic identity. \emph{Pacific J. Math.} {\bf 10}  (1960), 731-742.

\bibitem{BC} J. Baez and A.S. Crans, Higher-Dimensional Algebra VI: Lie $2$-Algebras. \emph{Theory and Appl. Categ}, {\bf 12}  (2004), 492-528.

\bibitem{BHC}
J. Baez, A. Hoffnung and C. Rogers, Categorified symplectic geometry and the
classical string. \emph{Comm. Math. Phys.} {\bf 293}  (2010), 701-725.


\bibitem{BR}J. Baez and C. Rogers, Categorified symplectic geometry and the string Lie 2-algebra. \emph{Homology, Homotopy Appl.}, {\bf 12} (2010), 221-236.

\bibitem{BSZ}
C. Bai, Y. Sheng and C. Zhu, Lie $2$-bialgebras. \emph{Comm. Math. Phys.} {\bf 320} (2013), 149-172.

\bibitem{Bax} R. J. Baxter, One-dimensional anisotropic Heisenberg chain. \emph{Ann. Physics} {\bf 70}  (1972), 323-337.

\bibitem{Pre-lie algebra in geometry} D. Burde, Left-symmetric algebras and pre-Lie algebras in geometry and physics, {\em Cent. Eur. J. Math.} {\bf 4} (2006), 323-357.

\bibitem{CP} V. Chari and A. Pressley, A Guide to Quantum Groups, Cambridge University Press, 1994.

\bibitem{CK}
A. Connes and D. Kreimer, { Renormalization in quantum field theory and the Riemann-Hilbert problem. I. The Hopf algebra structure of graphs and the main theorem.} {\em Comm. Math. Phys.} {\bf 210} (2000), 249-273.

\bibitem{Gub}
L. Guo,  An introduction to Rota-Baxter algebra. Surveys of Modern Mathematics, 4. International Press, Somerville, MA; Higher Education Press, Beijing, 2012. xii+226 pp.

\bibitem{JS}J. Jiang and Y. Sheng, Representations and cohomologies of relative Rota-Baxter Lie algebras and applications, arXiv:2108.08294.

\bibitem{Rota}    G.-C. Rota, Baxter algebras and combinatorial identities I, II. \emph{Bull. Amer. Math. Soc.} {\bf 75}  (1969), 325-329, 330-334.

\bibitem{LS}
T. Lada and J. Stasheff, Introduction to sh Lie algebras for
physicists. \emph{Internat. J. Theoret. Phys.} {\bf 32}  (1993), 1087-1103.


\bibitem{LST}
A. Lazarev, Y. Sheng and R. Tang, Deformations and homotopy theory of relative Rota-Baxter Lie algebras. \emph{Comm. Math. Phys.} {\bf 383}  (2021), 595-631.

\bibitem{Roy}
D. Roytenberg, Courant algebroids and strongly homotopy Lie algebras. \emph{Lett. Math. Phys.}, {\bf46}  (1998), 81-93.

\bibitem{Semonov-Tian-Shansky}
M. Semonov-Tian-Shansky, What is a classical R-matrix? \emph{Funct. Anal. Appl.} {\bf 17}(1983), 259-272 .

\bibitem{Sh} Y. Sheng, Categorification of pre-Lie Algebras and solutions of 2-graded Classical Yang-Baxter equations, \emph{Theor. Appl. Categ.}{\bf 34} (2019), No.11, 269-294.

\bibitem{TBGS}
R. Tang, C. Bai, L. Guo and Y. Sheng, Deformations and their controlling cohomologies of $\huaO$-operators. \emph{Commun. Math. Phys.} {\bf 368}  (2019), 665-700.

\bibitem{Yang} C. N. Yang, Some exact results for the many-body problem in one dimension with repulsive deltafunction
interaction. \emph{Phys. Rev. Lett.} {\bf19}  (1967), 1312-1315.


\end{thebibliography}
\end{document}